\documentclass[12pt]{amsart}
\usepackage[utf8]{inputenc}

\usepackage{mathrsfs} 
\usepackage{amssymb}
\usepackage{bm}
\usepackage{graphicx}
\usepackage[centertags]{amsmath}
\usepackage{amsfonts}
\usepackage{amsthm}
\usepackage{enumerate}
\usepackage{tocvsec2}
\usepackage{xcolor}
\usepackage[margin=.8in]{geometry}

\newtheorem{theorem}{Theorem}[section]
\newtheorem*{theorem*}{Theorem}

\newtheorem{corollary}[theorem]{Corollary}
\newtheorem{lemma}[theorem]{Lemma}
\newtheorem{rem}[theorem]{Remark}

\newtheorem{proposition}[theorem]{Proposition}

\newtheorem*{fact*}{Fact}

\theoremstyle{definition}
\newtheorem{definition}[theorem]{Definition}

\newcommand{\ee}{\varepsilon}
\newcommand{\nn}{\mathbb{N}}
\newcommand{\rr}{\mathbb{R}}

\newcommand{\supp}{\text{supp}}

\begin{document}

\title[Equivalence of block sequences]{Equivalence of block sequences in Schreier spaces and their duals}
\author{R.M.~Causey}
\address{R.M~Causey}
\email{rmcausey1701@gmail.com}

\author{A.~Pelczar-Barwacz}
\address{Jagiellonian University, Faculty of Mathematics and Computer Science, Institute of Mathematics, \L ojasiewicza 6, 30-348 Krak\'{o}w, Poland}
\email{anna.pelczar@uj.edu.pl}


\thanks{The research of the second author was supported by the grant of the National Science Centre (NCN), Poland, no. UMO-2020/39/B/ST1/01042}

\begin{abstract} We prove that any normalized block sequence in a Schreier space $X_\xi$, of arbitrary order $\xi<\omega_1$, admits a subsequence equivalent to a subsequence of the canonical basis of some Schreier space. The analogous result is proved for dual spaces to Schreier spaces. Basing on these results, we examine the structure of strictly singular operators on Schreier spaces and show that there are $2^\mathfrak{c}$ many closed operator ideals on a Schreier space of any order, its dual and bidual space.  
\end{abstract}

\maketitle

\section{Introduction}
The Schreier space was in introduced in \cite{Schr} providing an example of a basic sequence that is weakly null but has no Ces\`aro summing subsequence, which answered the question of Banach and Saks. Its definition relies on the Schreier family $\mathcal{S}_1$ of finite subsets of integers, given by the formula $\mathcal{S}_1=\{A\subset \nn: |A|\leq \min A\}$. The hierarchy of Schreier families $\mathcal{S}_\xi$, $\xi<\omega_1$, introduced in \cite{AA}, lies at the heart of the theory of mixed Tsirelson spaces and spaces built on their basis, cf. \cite{ADKM,AT}. Combinatorial spaces built with the use of families $\mathcal{S}_\xi$, $\xi<\omega_1$, are called Schreier spaces $X_\xi$, $\xi<\omega_1$. This model class of combinatorial spaces is studied widely, in particular with respect to  isometric geometric properties \cite{ABC,BC,BFT,CGS} and the structure of the lattice of closed operator ideals \cite{BKL,MP}, or used as reference spaces for properties of general Banach spaces \cite{C1}. Schreier spaces and other combinatorial spaces also form the basis for further constructions of spaces and a source of counterexamples to various conjectures on the structure of Banach spaces, cf. \cite{BL,CG}. 

The paper \cite{GL} contains a thorough study of the structure of block sequences in finite order Schreier spaces $X_n$, $n\in\nn$, providing the framework for the further study of bounded operators between subspaces spanned by subsequences of the basis and the rich family of pairwise non-comparable complemented subspaces. The results of \cite{GL} form crucial tools in the research in the structure of the operator algebra in \cite{BKL,MP}, in this case limited to the case of Schreier spaces of finite order. 

Following this approach, we characterize the family of basic sequences that is "sequence-satura\-ting" block sequences in a Schreier space of arbitrary order and its dual space. We show that any normalized block sequence in $X_\xi$, for any $\xi<\omega_1$, contains a subsequence equivalent to a subsequence of the canonical basis $(e_i^\varrho)_{i=1}^\infty$ of some Schreier space $X_\varrho$ (with an explicit recipe for $\varrho$), spanning a complemented subspace of $X_\xi$ and paired with a block sequence in $X_\xi^*$ with analogous properties. Moreover, we identify (in terms of $\xi$) those $\varrho<\omega_1$, for which such representation of subsequences of $(e_i^\varrho)_{i=1}^\infty$ is realized in $X_\xi$, as members of the set $R(\xi)=\{\xi-\iota: \iota\leqslant\xi\}$.  Note that the (finite) set $R(\xi)$ can be easily recovered from the Cantor normal form of the ordinal $\xi$. An analogous result is proved for dual spaces  $X_\xi^*$, $\xi<\omega_1$. 

The sequence pairing mentioned above takes the following form. Given a Banach space $X$ and  $\varrho<\omega_1$, we call a pair of normalized basic sequences  $((x_i)_{i=1}^\infty,(x_i^*)_{i=1}^\infty)\subset X\times X^*$ a $\varrho$-Schreier pair, provided that for some $J\in[\nn]$ we have the following.
\begin{enumerate}[$(i)$]\item $x^*_i(x_j)=0$ for all $i,j\in\nn$ with $i\neq j$, 
\item $\text{Re\ }x^*_i(x_i)=|x^*_i|(|x_i|)$ for all $i\in\nn$, 
\item $\inf_i \text{Re\ } x^*_i(x_i)>0$, 
\item $(x_i)_{i=1}^\infty$ is equivalent to a subsequence $ (e^{\varrho}_i)_{i\in J}$ of the canonical basis of $X_\varrho$, 
\item $(x^*_i)_{i=1}^\infty$ is equivalent to a subsequence  $(f^{\varrho}_i)_{i\in J}$ of the canonical basis of $X_\varrho^*$, 
\item $\sup_n \|\sum_{i=1}^n x_i\otimes x^*_i\|<\infty$. 
\end{enumerate}
It follows easily that sequences satisfying the above conditions span complemented subspaces in $X_\xi$ and $X_\xi^*$, respectively. Having this notion, the main results of the present paper, on equivalence of block sequences, can be then summarized as follows (see Theorems \ref{final1} and \ref{final2} for Theorem \ref{main1}, and Corollaries \ref{moveup1}, \ref{moveup2} and Theorem \ref{existence} for Theorem \ref{main2}). 
\begin{theorem}
Fix $\xi<\omega_1$. 
\begin{enumerate}
\item Any normalized block sequence in $X_\xi$ contains a subsequence $(x_i)_{i=1}^\infty$ such that for some  normalized block sequence $(x_i^*)_{i=1}^\infty\subset X_\xi^*$ and $\varrho\in R(\xi)$, $((x_i)_{i=1}^\infty,(x_i^*)_{i=1}^\infty)$ is a $\varrho$-Schreier pair.
\item Any normalized block sequence in $X_\xi^*$ contains a subsequence $(x_i^*)_{i=1}^\infty$ such that for some  normalized block sequence $(x_i)_{i=1}^\infty\subset X_\xi$ and $\varrho\in R(\xi)$, $((x_i)_{i=1}^\infty,(x_i^*)_{i=1}^\infty)$ is a $\varrho$-Schreier pair.
\end{enumerate}
\label{main1}
\end{theorem}

\begin{theorem}
Fix $\xi,\varrho <\omega_1$. The following are equivalent.
\begin{enumerate}
\item $\varrho\in R(\xi)$.
\item There is a normalized block sequence $(x_i)_{i=1}^\infty\subset X_\xi$ equivalent to a subsequence of the canonical basis of $X_\varrho$.
\item There is a normalized block sequence $(x_i^*)_{i=1}^\infty\subset X_\xi^*$ equivalent to a subsequence of the canonical basis of $X_\varrho^*$.
\item There is a $\varrho$-Schreier pair $((x_i)_{i=1}^\infty,(x_i^*)_{i=1}^\infty)\subset X_\xi\times X_\xi^*$. 
\end{enumerate}
If any of the above conditions is satisfied, the sequences $(x_i)_{i=1}^\infty$ and $(x_i^*)_{i=1}^\infty$ in (4) can be chosen so that each $x_i$,  $i\in\nn$, is a convex combination of elements of the canonical basis of $X_\xi$, and each $x_i^*$,  $i\in\nn$, is a finite sum of elements of the canonical basis of $X_\xi^*$. 
\label{main2}
\end{theorem}
Concerning Theorem \ref{main1}, on the level of picking subsequences of normalized block sequences in $X_\xi$ or $X_\xi^*$, that are equivalent to some subsequence of the basis of some $X_\varrho$ or $X_\varrho^*$, the case of finite order Schreier spaces $X_n$, $n\in\nn$, was solved in \cite[Prop. 4.3]{GL}, the case of $X_1^*$ was proved in \cite{L}. The proof of  Theorems \ref{main1} and \ref{main2} requires technical combinatorial tools on the Schreier families, in particular the equality of any Schreier family of arbitrary order and its modified version (Lemma \ref{mod1}). This result was proved first for finite order in \cite{ADKM} and shown for arbitrary order in \cite{Schl2}, we present an alternative proof. Identification of the order $\varrho$ of a Schreier space  $X_\varrho$ to be represented on a given block sequence in $X_\xi$ is given in the spirit of \cite{GL}, whereas in the case of dual Schreier spaces is more technical and requires additional notions of $\xi$-flat and $\xi$-approachable sequences of functionals (see Section 5.).

The structure of block sequences in a Banach space is closely related to the complexity of the class of strictly singular operators. Recall that an operator $T:X\to Y$ is strictly singular, if for any $\varepsilon>0$ in any infinite dimensional subspace of $X$ there is a normalized vector $x$ with $\|Tx\|<\varepsilon$. The appropriate framework in our case is described by classes of $\varrho$-strictly singular operators, $\varrho<\omega_1$, introduced in \cite{ADST}, with strict singularity measured by means of Schreier families $\mathcal{S}_\varrho$.  We apply Theorems \ref{main1} and \ref{main2} to prove that for any $\xi<\omega_1$, the classes of $(\varrho+1)$-strictly singular operators on $X_\xi$, $\varrho\in R(\xi)$, are pairwise distinct (closed) operator ideals, exhausting the ideal of strictly singular operators (Theorem \ref{structure of ss operators}). Moreover, any product of $|R(\xi)|$ many strictly singular operators is compact.  

Finally, we apply the combinatorial results on Schreier families to prove that there are $2^\mathfrak{c}$ many pairwise distinct closed operator ideals on Schreier spaces of arbitrary order, their dual and bidual spaces (Theorem \ref{cardinality}, concerning Schreier spaces see \cite{Schl}), in particular Schreier spaces and their duals are on the list of separable spaces with the maximal cardinality of the lattice of closed operator ideals, which includes spaces $L_p(0,1)$, $1<p<\infty$, $\ell_p\oplus\ell_q$, $1\leq p<q\leq\infty$, and $\ell_p\oplus c_0$, $1\leq p<\infty$, products of reflexive $\ell_p$ spaces and convexified Tsirelson spaces and their duals \cite{JS, FSZ}. The reasoning is based on generalizations of certain techniques of \cite{GL,MP}, whose proof relies on  combinatorial  Lemmas \ref{mod1}  and \ref{weaksumming}. 

The paper is organized as follows. In Section 2. we recall the notion of  Schreier families $\mathcal{S}_\xi$, $\xi<\omega_1$ along with their basic properties and prepare the combinatorial tools for the next part. In Section 3. we recall the definition and basic properties of Schreier spaces, their dual and bidual spaces, and introduce the notion of Schreier pair of basic sequences. Section 4. is devoted to the crucial tool we shall use in our reasoning, i.e. the hierarchy of repeated averages, with the central  Lemma \ref{weaksumming} and isometric representation of $c_0$ in any Schreier space (Corollary \ref{isometricc0}). In Section 5. we introduce the notion of  $\xi$-flat and $\xi$-approachable sequences of functionals on Schreier spaces, providing tools for identifying the order $\varrho\in R(\xi)$ of a dual Schreier space to be represented on a given  block sequence in the dual Schreier space.  Section 6. contains the main results characterizing basic sequences "sequence-saturating" block sequences in Schreier spaces and their dual spaces.  Section 7. is devoted to classes of strictly singular operators on Schreier spaces, Section 8. concerns the cardinality of the lattice of closed operator ideals on Schreier spaces of arbitrary order, their dual and bidual spaces. 

\section{Schreier families}

In this section we recall the basic notation, the hierarchy of Schreier families introduced in \cite{AA} and their standard properties.  We recall the fact that Schreier families and their modified versions coincide (\cite{Schl2}) with an alternative proof. We end the section with a number of combinatorial observations providing relations between Schreier families of different order we shall use in the sequel.

Throughout, $\mathbb{K}$ will denote the scalar field, which is either $\mathbb{R}$ or $\mathbb{C}$.  For an ordinal $\xi$, $\xi=[0,\xi)$ and $\xi+=\xi+1=[0,\xi]$.   

Subsets of $\nn$ are identified with strictly increasing sequences of natural number by identifying a set with the sequence obtained by listing the members of the set in increasing order. Therefore we use both set and sequence notation interchangeably. In particular, for $n\in\nn$, we let $(n)$ denote the singleton containing $n$. 

For any ordinal $\xi$ of the form $\xi=\omega^{\beta_1}+\dots+\omega^{\beta_k}$, with $k\in\nn$ and ordinals $\beta_1\geqslant\dots\geqslant\beta_k\geqslant 0$ we write \[I(\xi):=\{0\}\cup\{\omega^{\beta_1}+\dots+\omega^{\beta_j}: 1\leqslant j\leqslant k\}, \ R(\xi):=\{\xi-\iota: \iota\in I(\xi)\} =\{\xi-\iota: \iota\leqslant\xi\}.\]

We let $2^\nn$, $[\nn]^{<\omega}$, and $[\nn]$ denote the sets of all subsets of $\nn$, all finite subsets of $\nn$, and all infinite subsets of $\nn$, respectively.  For $M\in[\nn]$, we let $2^M, [M]^{<\omega}$, and $[M]$ denote the sets of all subsets of $M$, all finite subsets of $M$, and all infinite subsets of $M$, respectively. 

We establish the conventions that $\max \varnothing=0$ and $\min \varnothing=\infty$. We use $E<F$ to denote the relation that $\max E<\min F$. We use $\smallfrown$ to denote concatenation of sequences.   For $E\in 2^\nn$ and $i\in \{j\in\nn:j\leqslant |E|\}$, we let $E(i)$ denote the $i^{th}$ smallest member of $E$. That is, if $E\in[\nn]^{<\omega}$, $E=\{E(1), \ldots, E(|E|)\}$ with $E(1)<\ldots <E(|E|)$, and if $M\in[\nn]$, $M=\{M(1), M(2), \ldots\}$ with $M(1)<M(2)<\ldots$.  Given any $E\subset\nn$ and $M\in[\nn]$, we let $M(E)=\{M(i):i\in E\}$. We let $E\preceq F$ denote the relation that $E$ is an initial segment of $F$, which is equivalent to the condition that $|E|\leqslant |F|$ and $E(i)=F(i)$ for all $1\leqslant i\leqslant |F|$.    Given $\mathcal{F}\subset [\nn]^{<\omega}$ and $M\in[\nn]$, we let 
\[\mathcal{F}(M)=\{M(E): E\in\mathcal{F}\}.\]
 For any $\mathcal{F}\subset [\nn]^{<\omega}$ we denote by $MAX (\mathcal{F})$ the family of elements of $\mathcal{F}$ that are maximal with respect to inclusion and we let \[\mathcal{F}'=\mathcal{F}\setminus MAX(\mathcal{F}).\] 
We  define the \emph{Schreier families} $\mathcal{S}_\xi$, $\xi<\omega_1$ (\cite{AA}).  We define \[\mathcal{S}_0=\{\varnothing\}\cup \{(n):n\in\nn\}.\]   If $\mathcal{S}_\xi$ has been defined, we let \[\mathcal{S}_{\xi+1}= \{\varnothing\}\cup \Bigl\{\bigcup_{i=1}^k E_i: k\leqslant E_1<\ldots <E_k, \varnothing\neq E_i\in \mathcal{S}_\xi\Bigr\}.\] If $\xi<\omega_1$ is a limit ordinal, we fix a sequence $\xi_n\uparrow \xi$ such that $\mathcal{S}_{\xi_n+1}\subset \mathcal{S}_{\xi_{n+1}}$ for all $n\in\nn$ (cf. \cite[Section 3.2]{C2}) and define \[\mathcal{S}_\xi=\{\varnothing\}\cup \{E\in[\nn]^{<\omega}\setminus\{\varnothing\}:\exists n\leqslant E\in \mathcal{S}_{\xi_n+1}\}= \{\varnothing\}\cup \{E\in[\nn]^{<\omega}\setminus\{\varnothing\}: E\in \mathcal{S}_{\xi_{\min E}+1}\}.\] 
For each $\xi<\omega_1$ and $n\in\nn$, we define \[\mathcal{A}_n[\mathcal{S}_\xi]=\{\varnothing\}\cup\Bigl\{\bigcup_{i=1}^k E_i:k\leqslant n, E_1<\ldots <E_k, \varnothing\neq E_i\in \mathcal{S}_\xi\Bigr\}.\] 
We next define the \emph{modified Schreier families} $\mathcal{S}_\xi^M$, $\xi<\omega_1$.  We define \[\mathcal{S}_0^M=\{\varnothing\}\cup \{(n):n\in\nn\}.\]   If $\mathcal{S}_\xi^M$ has been defined, \[\mathcal{S}_{\xi+1}^M= \{\varnothing\}\cup \Bigl\{\bigcup_{i=1}^k E_i: k\leqslant E_i \in \mathcal{S}_\xi^M\setminus \{\varnothing\}, E_1,\dots,E_k\text{\ pairwise disjoint}\Bigr\}.\]  If $\xi<\omega_1$ is a limit ordinal, we fix the same sequence $\xi_n\uparrow \xi$ used in the definition of $\mathcal{S}_\xi$ an define  \[\mathcal{S}_\xi^M=\{\varnothing\}\cup \{E\in[\nn]^{<\omega}\setminus\{\varnothing\}:\exists n\leqslant E\in \mathcal{S}_{\xi_n+1}^M\}.\]

\begin{definition} For $(l_i)_{i=1}^n,(m_i)_{i=1}^n\in [\nn]^{<\omega}$, we say $(m_i)_{i=1}^n$ is a \emph{spread} of $(l_i)_{i=1}^n$ if $l_i\leqslant m_i$ for all $1\leqslant i\leqslant n$.  We say $\mathcal{F}\subset [\nn]^{<\omega}$ is \emph{spreading} if it contains all spreads of its members.  We say $\mathcal{F}\subset [\nn]^{<\omega}$ is \emph{hereditary} if it contains all subsets of its members. We identify $E$ with its indicator function $1_E\in \{0,1\}^\nn$. We endow $\{0,1\}^\nn$ with the product topology and identify $2^\nn$ with the topology induced by the identification $E\leftrightarrow 1_E$.  We say $\mathcal{F}\subset [\nn]^{<\omega}$ is \emph{compact} if it is compact with respect to this topology. We say $\mathcal{F}\subset [\nn]^{<\omega}$ is \text{Regular} if it is compact, spreading, and hereditary.

\end{definition}

 We recall now the following facts, $(i)$-$(iv)$ of which are well known, cf. \cite{AMT,AT}

\begin{proposition} [Properties of Schreier families] Fix $\xi<\omega_1$. 

\begin{enumerate}[(i)]\item For all $\xi<\omega_1$ and $n\in\nn$, $\mathcal{A}_n[\mathcal{S}_\xi]$ is regular. 

\item For any $\xi<\omega_1$ and $n\in\nn$, a member $E$ of $\mathcal{A}_n[\mathcal{S}_\xi]$ is maximal in $(\mathcal{A}_n[\mathcal{S}_\xi],\subset)$ if and only if it is maximal in $(\mathcal{A}_n[\mathcal{S}_\xi],\preceq)$. 

\item For $E\in [\nn]^{<\omega}$, there exists $E<m\in\nn$ such that $E\smallfrown(m)\in \mathcal{S}_\xi$ if and only if for all $E<m\in\nn$, $E\smallfrown(m)\in\mathcal{S}_\xi$. 

\item For any $\beta\leqslant \xi<\omega_1$, there exists $m=m(\beta,\xi)$ such that $\mathcal{S}_\beta\cap [[m,\infty)]^{<\omega}\subset \mathcal{S}_\xi$. 

\item $\mathcal{S}_1\subset \mathcal{S}_\xi$ for all $1\leqslant \xi<\omega_1$. 

\item For any $E\in MAX(\mathcal{S}_\xi)$, $\max E\geqslant 2\min E-1$. 

\item If $E_1<E_2<\ldots$ are successive, maximal members of $\mathcal{S}_\xi$ for some $1\leqslant \xi<\omega_1$, then for all $i,j\in \nn$, \[\frac{\max E_i}{\min E_{i+j}} \leqslant 2^{1-j}.\]

\end{enumerate}
\label{schreierfamilyfacts}
\end{proposition}

\begin{proof}
$(v)$ We prove by induction on $\xi\in [1,\omega_1)$ that $\mathcal{S}_1\subset \mathcal{S}_\xi$.   Of course, the $\xi=1$ case is trivial. If $\mathcal{S}_1\subset \mathcal{S}_\xi$, then $\mathcal{S}_1\subset \mathcal{S}_\xi\subset \mathcal{S}_{\xi+1}$. Last, fix $\xi<\omega_1$ and fix $\xi_n\uparrow \xi$.  For $\varnothing\neq E\in \mathcal{S}_1$, since $1\leqslant \xi_1+1$, $1\leqslant E\in \mathcal{S}_{\xi_1+1}$, so $E\in \mathcal{S}_\xi$.  

$(vi)$  Fix $1\leqslant \xi<\omega_1$ and $E\in MAX(\mathcal{S}_\xi)$. Since $\mathcal{S}_1\subset \mathcal{S}_\xi$ (Prop. \ref{schreierfamilyfacts} (v)), there is an initial segment $F$ of $E$ such that $F\in MAX(\mathcal{S}_1)$. But this means that $|F|=\min F$, and \[\max E\geqslant \max F\geqslant \min F-1+|F|=2\min F-1=2\min E-1.\]

$(vii)$ Let  $E_1<E_2<\ldots$ be successive, maximal members of $\mathcal{S}_\xi$ for some $1\leqslant \xi<\omega_1$. For any $k\in\nn$, \[\min E_{k+1}\geqslant 1+\max E_k \geqslant 1+2\min E_k-1=2\min E_k,\] so \[\min E_{i+j}\geqslant 2\min E_{i+j-1}\geqslant 2^2 \min E_{i+j-2}\geqslant \ldots \geqslant  2^{j-1}\min E_{i+1}>2^{j-1}\max E_i.\] 

\end{proof}

The following was shown first for $\xi<\omega$  in \cite{ADKM}  and proved for any $\xi<\omega_1$  in \cite{Schl2}.  Using our specific choices of sequences $\xi_n\uparrow \xi$ for countable limit ordinals $\xi$, we provide an alternative proof for any $\xi<\omega_1$.

\begin{lemma}  For each $\xi<\omega_1$, $\mathcal{S}_\xi=\mathcal{S}_\xi^M$. 
\label{mod1}
\end{lemma}

\begin{proof} 
For each $\xi<\omega_1$ and $k\in\nn$, let $P_\xi(k)$ be the proposition that for any non-empty, pairwise disjoint $E_1, \ldots, E_k$ in $\mathcal{S}_\xi$ with $\min E_1<\ldots <\min E_k$, there exist $E_1'<\ldots <E_k'$ in $\mathcal{S}_\xi$ such that $\cup_{i=1}^k E_i=\cup_{i=1}^k E_i'$ and $\min E_i\leqslant \min E_i'$ for all $1\leqslant i\leqslant k$.  Let $P_\xi$ be the proposition that $P_\xi(k)$ holds for all $k\in \nn$. Let $Q_\xi$ be the proposition that $\mathcal{S}_\xi=\mathcal{S}_\xi^M$.  Let $R_\xi$ be the proposition that $P_\xi$ and $Q_\xi$ hold.

We shall prove that for each $\xi<\omega_1$, $R_\xi$ holds.  It is quite clear that for any $\xi<\omega_1$, if we have $P_\zeta$ for all $\zeta\leqslant \xi$, then $Q_\xi$ follows.  Therefore once we have that $P_\zeta$ holds for all $\zeta\leqslant \xi$, we will have $R_\xi$.    We work by induction, for which it is sufficient to prove that for each $\xi<\omega_1$,  if $R_\zeta$ holds for all $\zeta<\xi$, then $P_\xi$ holds. 

Obviously $P_0$ holds.   

Assume that for some $\xi<\omega_1$, $P_\xi$ holds. We prove that $P_{\xi+1}$ holds. The proof is identical to the proof from \cite{ADKM} to establish the cases $\xi<\omega$. We repeat the proof here for completeness.  Fix non-empty, pairwise disjoint sets $E_1, \ldots, E_k\in \mathcal{S}_{\xi+1}$ with $\min E_1<\ldots < \min E_k$.   For each $1\leqslant i\leqslant k$, we can write $E_i=\cup_{j=1}^{m_i}E_{i,j}$ for some $E_{i,1}<\ldots < E_{i,m_i}$, $\varnothing\neq E_{i,j}\in \mathcal{S}_\xi$, $m_i\leqslant \min E_i$.   Let $s_i=\sum_{j=1}^i m_j$ for $i=0, \ldots, k$.   Let $(F_j)_{j=1}^{s_k}$ be an enumeration of $(E_{i,j})_{i=1,j=1}^{k,m_i}$ such that $\min F_1<\ldots < \min F_{s_k}$.   By the inductive hypothesis, there exist $F_1'<\ldots <F_{s_k}'$ in $\mathcal{S}_\xi$ such that $\min F_j\leqslant \min F_j'$ for all $1\leqslant j\leqslant s_k$ and $\cup_{i=1}^k E_i=\cup_{j=1}^{s_k}F_j=\cup_{j=1}^{s_k}F_j$.    Define \[E_i'=\bigcup_{j=s_{i-1}+1}^{s_i}F_j'\] for all $1\leqslant i\leqslant k$.   It is easy to see that \[m_i\leqslant \min E_i \leqslant \min F_{s_{i-1}+1} \leqslant \min F_{s_{i-1}+1}'=\min E_i' \]  for all $1\leqslant i\leqslant m$.  Since $\min E_i'\geqslant m_i$ and $E_i'$ is the union of $m_i$ successive members of $\mathcal{S}_\xi$, $E_i'\in \mathcal{S}_{\xi+1}$ for all $1\leqslant i\leqslant k$. This finishes the successor case.

Last, let $\xi<\omega_1$ be a limit ordinal and assume $R_\zeta$ holds for all $\zeta<\xi$.     Let $\xi_n\uparrow \xi$ be the sequence used to define $\mathcal{S}_\xi$.   We prove $P_\xi$ by proving $P_\xi(k)$ by induction on $k$.  The case $k=1$ is trivial.    Assume that for some $1<k$, the result holds for $k-1$ and fix non-empty, pairwise disjoint $E_1, \ldots, E_k$ in $\mathcal{S}_\xi$ with $\min E_1<\ldots < \min E_k$.   If $\max E_1<\min E_2$, then $E_1<\cup_{i=2}^{k-1}E_i$. We let $E_1'=E_1$ and apply $P_\xi(k-1)$ to $E_2, \ldots, E_k$ to obtain $E_2', \ldots, E_k'$ to finish the proof. Assume $\max E_1>\min E_2$.  Note that $E_i\in \mathcal{S}_{\xi_{\min E_i}+1}$ for all $1\leqslant i\leqslant k$. 

 Fix $2\leqslant l\leqslant k$ such that $\min E_i\leqslant \max E_1$ for all $2\leqslant i\leqslant l$ and $\min E_i>\max E_1$ for all $l<i\leqslant k$.   Define \[F_i=\left\{\begin{array}{ll} (E_1\cup (\min E_2))\setminus (\max E_1) & : i=1 \\ (E_i\cup (\min E_{i+1}))\setminus (\min E_i) & : 2\leqslant i < l \\ (E_i\cup (\max E_1))\setminus (\min E_i) & : i=l \\  E_i & : l < i \leqslant k. \end{array}\right.\]  Note that 

\begin{enumerate}[(a)]\item $F_1, \ldots, F_k$ are pairwise disjoint, 

\item $\min F_1<\ldots<\min F_k$, 

\item $\varnothing\neq F_i$ for all $1\leqslant i\leqslant k$, 

\item $F_i$ is a spread of $E_i$ for all $2\leqslant i\leqslant k$, 

\item $F_i\in \mathcal{S}_{\xi_{\min E_i}+1}$ for all $2\leqslant i\leqslant k$, 

\item $\min E_i\leqslant \min F_i$ for all $1\leqslant i\leqslant k$.

\end{enumerate}

Define $E_1'=F_1\cap [1, \max E_2]<\cup_{i=2}^k F_i$ and note that $\min E_1=\min E_1'$ and $E_1'\in \mathcal{S}_{\xi_{\min E_1}+1}$, as noted in Proposition \ref{schreierfamilyfacts}. So $E_1'\in \mathcal{S}_\xi$.

Let $G_1=F_2\cup (F_1\setminus E_1')$ and let $G_i=F_{i+1}$ for $1<i<k$.  Note that $\min G_1<\ldots <G_{k-1}$, $\min G_i\geqslant \min E_{i+1}$ for all $1\leqslant i<k$, and $G_i\in \mathcal{S}_\xi$ for $1<i<k$.    In order to apply $P_\xi(k-1)$ to $G_1, \ldots, G_{k-1}$, we only need to show that $G_1\in \mathcal{S}_\xi$.   Since $E_2\in \mathcal{S}_{\xi_{\min E_2}+1}$, we can write $E_2=\cup_{j=1}^p H_j$ for some successive, and therefore pairwise disjoint, $H_1, \ldots, H_p\in \mathcal{S}_{\xi_{\min E_2}}$ and $p\leqslant \min E_2$.    Since $F_2$ is a spread of $E_2$, there exist spreads $I_1, \ldots, I_p$ of $H_1, \ldots, H_p$, respectively, such that $F_2=\cup_{j=1}^p I_j$.   Since $E_1\in \mathcal{S}_{\xi_{\min E_1}+1}\subset \mathcal{S}_{\xi_{\min E_2}}$ and \[\min \{\min F_2, \min [F_1\setminus E_1']\}>\min E_2,\] $I_1, \ldots, I_p, [F_1\setminus E_1']$ are pairwise disjoint subsets of $\mathcal{S}_{\xi_{\min E_2}}\cap [[p+1,\infty)]^{<\omega}$, and therefore $G_1=\bigl(\cup_{j=1}^p I_j\bigr)\cup [F_1\setminus E_1']\in \mathcal{S}_{\xi_{\min E_2}+1}^M=\mathcal{S}_{\xi_{\min E_2}+1}$. Here we have used $R_{\xi_{\min E_2}+1}$.   Since $\min E_2\leqslant G_1$, $G_1\in \mathcal{S}_\xi$.   

 Therefore we can apply $P_\xi(k-1)$ to $G_1, \ldots, G_{k-1}$ to get $E_2'<\ldots <E_k'$ in $\mathcal{S}_\xi$ such that $\min E_{i+1}'\geqslant \min G_i\geqslant \min F_{i+1}\geqslant \min E_i$ for all $2\leqslant i\leqslant k$ and $\cup_{i=1}^{k-1}G_i=\cup_{i=2}^k E_i'$. Of course, $E_1'<E_2'$, so $E_1',  \ldots, E_k'$ satisfy the conclusions. 
 
 \end{proof}

In the sequel we shall need the following observations relating Schreier families of different order, making use of our specific choice of sequences of ordinals in the definition of Schreier families (the limit case).

\begin{proposition} 
\begin{enumerate}[(i)]\item Suppose that for some ordinals $ \beta, \delta, \varrho$, with $\beta<\omega^\delta\leqslant \varrho$. Then $\beta+\varrho=\varrho$.   

\item For any $\gamma,\varrho<\omega_1$ with $\varrho\geqslant 1$, any $n\in\nn$,  any $E_1<E_2<\ldots$, $i\leqslant  E_i\in  \mathcal{A}_n[\mathcal{S}_\gamma]$,  and any $K\in[\nn]$, there exists $L\in[K]$ such that \[\Bigl\{\bigcup_{j\in F}E_{L(j)}: F\in \mathcal{S}_\varrho\Bigr\}\subset \mathcal{S}_{\gamma+\varrho}.\] Moreover, the inclusion \[\Bigl\{\bigcup_{j\in F}E_{M(j)}:F\in \mathcal{S}_\varrho\Bigr\}\subset \mathcal{S}_{\gamma+\varrho}\] holds for any $M\in[L]$.   

\item Fix ordinals $\iota, \delta, \varrho$ with $\varrho\geqslant \omega^\delta$, a sequence of ordinals $(\beta_i)_{i=1}^\infty\subset \omega^\delta$, $N\in[\nn]$.  Suppose that $E_{i,j}$, $i,j\in\nn$ are such that for all $i\in\nn$, $E_{i,1}<E_{i,2}<\ldots$, and for all $i,j\in\nn$, $i\leqslant  E_{i,j}\in \mathcal{A}_{N(i)}[\mathcal{S}_{\iota+\beta_i}]$.  Then for any $K\in [\nn]$, there exists $L\in[\nn]$ such that for all $M\in[L]$, \[\bigcup_{i=1}^\infty \Bigl\{\bigcup_{j\in F}E_{i,M(j)}:F\in \mathcal{S}_\varrho\cap [[i,\infty)]^{<\omega}\Bigr\}\subset \mathcal{S}_{\iota+\varrho}.  \]

\item For any $\gamma,\delta<\omega_1$, there exists $M\in[\nn]$ such that for any $E\in \mathcal{S}_{\gamma+\delta}$, there exist $F_1<\ldots <F_k$, $\varnothing\neq F_j\in \mathcal{S}_\gamma$ such that $E=\cup_{j=1}^k F_j$ and $M((\min F_j)_{j=1}^k)\in \mathcal{S}_\delta$. 

\item For any $\gamma,\delta<\omega_1$, there exists $N\in[\nn]$ such that for any $\varnothing\neq F\in \mathcal{S}_\delta$ and any successive sets $(E_i)_{i\in F}\subset \mathcal{S}_\gamma$ with $\min E_i\geqslant N(i)$ for all $i\in F$, then \[\bigcup_{i\in F}E_i\in \mathcal{S}_{\gamma+\delta}.\] 

\end{enumerate}
\label{combinatorics}
\end{proposition}

\begin{proof}$(i)$ It is a standard property of ordinals that $\beta+\varrho=\varrho$ provided $\beta\cdot\omega\leqslant \varrho$.

$(ii)$ Note that for $(ii)$, the only reason we must include the hypothesis that $i\leqslant E_i$ is that some of the $E_i$ may be empty. If each $E_i$ is non-empty, then the condition that $E_1<E_2<\ldots$ yields that $i\leqslant E_i$.  We note that the last sentence of $(ii)$ follows from the preceding parts, since $\mathcal{S}_\delta(M)\subset \mathcal{S}_\delta(L)$ for any $L\in[\nn]$ and $M\in[L]$.  Indeed, if $E\in \mathcal{S}_\delta(M)$, then $E=M(F)$ for some $F\in \mathcal{S}_\delta$.  If $j_1<j_2<\ldots $ are such that $M(i)=L(j_i)$ for all $i\in\nn$, then $E=M(F)=L(\{j_i:i\in F\})$. This set lies in $\mathcal{S}_\delta(L)$, since $\{j_i:i\in\nn\}$ is a spread of $F$, and is therefore a member of $\mathcal{S}_\delta$. 

We prove the first statement by induction on $\varrho$.  For the case $\varrho=1$, we choose any $L\in[K]$ such that $L(i)\geqslant in$.  Then for any $F\in \mathcal{S}_1$, $\cup_{j\in F}E_{L(j)}$ is the union of at most $n|F|\leqslant n \min F$ successive members of $\mathcal{S}_\gamma$ and $\min \cup_{j\in F}E_{L(j)}\geqslant L(\min F)\geqslant n\min F$.  Therefore $\cup_{j\in F}E_{L(j)}\in \mathcal{S}_{\gamma+1}$.   

Assume the result holds for some $\varrho$.  Fix $L\in[K]$ such that \[\Bigl\{\bigcup_{i\in F}E_{L(i)}: F\in \mathcal{S}_\varrho\Bigr\}\subset \mathcal{S}_{\gamma+\varrho}.\]  Then \[\Bigl\{\bigcup_{i\in F}E_{L(i)}: F\in \mathcal{S}_{\varrho+1}\Bigr\}\subset \mathcal{S}_{\gamma+\varrho+1}.\]  Indeed, for $F\in \mathcal{S}_{\varrho+1}$, write $F=\cup_{i=1}^m F_i$ with $F_1<\ldots <F_m$, $m\leqslant F_1$, and $\varnothing\neq F_i\in \mathcal{S}_\varrho$.   Let \[G_i=\bigcup_{j\in F_i}E_{L(j)}\] for each $1\leqslant i\leqslant m$.  By the inductive hypothesis, $G_i\in \mathcal{S}_{\gamma+\varrho}$ for each $1\leqslant i\leqslant m$.  Since $m\leqslant \min F_1\leqslant \min E_{\min F_1}=\min G_1$, $m\leqslant G_1<\ldots<G_m$, so \[\bigcup_{j\in F}E_{L(j)}=\bigcup_{i=1}^m G_i\in \mathcal{S}_{\gamma+\varrho+1}.\]   

Assume the result holds for all ordinals less than some limit ordinal $\varrho$.   Note that $\varrho$ and $\gamma+\varrho$ are both limit ordinals.   Fix $\varrho_i\uparrow \varrho$ and $\mu_i\uparrow\gamma+\varrho$ such that \[\mathcal{S}_\varrho=\bigcup_{i=1}^\infty \mathcal{S}_{\varrho_i+1}\cap[[i,\infty)]^{<\omega}\] and  \[\mathcal{S}_{\gamma+\varrho}=\bigcup_{i=1}^\infty \mathcal{S}_{\mu_i+1}\cap [[i,\infty)]^{<\omega}.\]  Fix $L_0\in[K]$ such that $\gamma+\varrho_i+1<\mu_{L_0(i)}+1$ for all $i\in\nn$ and such that $\mathcal{S}_{\gamma+\varrho_i+1}\cap [[L_0(i),\infty)]^{<\omega}$.   Recursively select $L_i\in [L_{i-1}]$ for $i\in\nn$ such that for all $i\in\nn$ and $M\in [L_i]$, \[\Bigl\{\bigcup_{j\in F}E_{L_i(j)}:F\in \mathcal{S}_{\varrho_i+1}\Bigr\}\subset \mathcal{S}_{\gamma+\varrho_i+1}.\]   Let $L(i)=L_i(i)$.   Fix $\varnothing\neq F\in \mathcal{S}_\varrho$ and let $i=\min F$.   Then $F\in \mathcal{S}_{\varrho_i+1}$.   Let $L'=\{L_i(1), \ldots, L_i(i-1), L(i), L(i+1), \ldots \}\in [L_i]$.    Therefore \[E:=\bigcup_{j\in F}E_{L(j)}=\bigcup_{j\in F}E_{L'(j)}\in \mathcal{S}_{\gamma+\varrho_i+1}.\]  Since \[\min E=\min E_{L(\min F)}=\min E_{L(i)}\geqslant L(i)\geqslant L_0(i),\] $E\in \mathcal{S}_{\gamma+\varrho_i+1}\cap [[L_0(i),\infty)]^{<\omega}\subset \mathcal{S}_{\mu_{L_0(i)}+1}$.   Again, since $\min E\geqslant L_0(i)$, \[E\in \mathcal{S}_{\mu_{L_0(i)}+1}\cap [[L_0(i),\infty)]^{<\omega}\subset \mathcal{S}_{\gamma+\varrho}.\] 

$(iii)$ By $(i)$, it follows that $\iota+\beta_i+\varrho=\iota+\varrho$ for all $i\in\nn$.    Let $L_0=K$. By repeated applications of $(ii)$, we can recursively select $L_i\in [L_{i-1}]$ such that for all $i\in\nn$ and $M\in[L_i]$, \[\Bigl\{\bigcup_{j\in F}E_{i,M(j)}:F\in \mathcal{S}_\varrho\Bigr\} \subset \mathcal{S}_{\iota+\beta_i+\varrho}=\mathcal{S}_{\iota+\varrho}.\]   Fix $M\in[L]$, $i\in \nn$, and $\varnothing \neq F\in \mathcal{S}_\varrho$ with $\min F\geqslant i$.   Let \[M'=\{L_i(1), \ldots, L_i(i-1), M(i), M(i+1), \ldots\}\in [L_i].\]  Then \[\bigcup_{j\in F}E_{i,M(j)}=\bigcup_{j\in F}E_{i,M'(j)}\in \mathcal{S}_{\iota+\varrho}.\] 

$(iv)$  We prove the claim by induction on $\delta$ for $\gamma$ fixed. For $\delta=0$ it suffices to take $M=\nn$.  

Assume that for $\delta<\omega_1$, $M\in[\nn]$ satisfies the conclusion. For  $E\in\mathcal{S}_{\gamma+\delta+1}$
we pick $l\leqslant G_1<\dots<G_l$ with  each  $G_j\in\mathcal{S}_{\gamma+\delta}$. 
By the inductive hypothesis, for any $j=1,\dots,l$ we have \[G_j:=\bigcup_{i\in A_j}F_i\] for some successive $(F_i)_{i\in A_j}\subset \mathcal{S}_{\gamma}$ with $M((\min F_i)_{i\in A_j})\in \mathcal{S}_{\delta}$. As $\min F_{\min A_j}=\min G_j$ for all $j=1,\dots,l$ we have $M((\min F_i)_{i\in \cup_jA_j})\in \mathcal{S}_{\delta+1}$, thus $M$ satisfies the conclusion also for $\delta+1$. 

Let $\delta$ be a limit ordinal and assume the result holds for all ordinals less than $\delta$.   Note that $\gamma+\delta$ is a limit ordinal.  Fix $\delta_n\uparrow\delta$ and $\mu_n\uparrow \gamma+\delta$ such that \[\mathcal{S}_\delta=\bigcup_{n=1}^\infty \mathcal{S}_{\delta_n+1}\cap [[n,\infty)]^{<\omega}\] and \[\mathcal{S}_{\gamma+\delta}=\bigcup_{n=1}^\infty \mathcal{S}_{\mu_n+1}\cap [[n,\infty)]^{<\omega}.\] 
Pick $I\in [\nn]$ such that for any $n\in\nn$, $\mu_n+1\leqslant\gamma+\delta_{I(n)}$ and 
\[\mathcal{S}_{\mu_n+1}\cap [[I(n),\infty)]^{<\omega}\subset \mathcal{S}_{\gamma+\delta_{I(n)}}.
\]
By the inductive assumption, for each $n\in\nn$ pick $M_n\in [\nn]$
satisfying the conclusion for $\delta_n$. Fix $M\in[\nn]$ such that for any $i\in\nn$ we have
\[M(i)\geqslant\max\{2,M_1(i), \dots,M_i(i), I(i)\}\]
and notice that $M([i,\infty))$ is a a spread of $M_n([i,\infty))$ for all $i\in\nn$. 
For a fixed $E\in\mathcal{S}_{\gamma+\delta}$, let $n=\min E$, then $E\in\mathcal{S}_{\mu_n+1}$. If $n=1$, then $E$ is a singleton and the conclusion follows, otherwise let $E'=E\cap [I(n),\infty)$. As $E'\in\mathcal{S}_{\gamma+\delta_{I(n)}}$, by the inductive assumption $E'=G_1\cup\dots\cup G_m$ for some $(G_i)_{i=1}^m\subset \mathcal{S}_\gamma$ with $M_{I(n)}((\min G_i)_{i=1}^m)\in\mathcal{S}_{\delta_{I(n)}}$. Let $E\setminus E'=\{n_1,\dots,n_l\}\subset \{n,\dots,I(n)\}$ and $(F_j)_{j=1}^k=(\{n_i\})_{i=1}^l\cup (G_i)_{i=1}^m$. Notice $M(\{n_1,\dots,n_l\})\in\mathcal{S}_1$  and so $I(n)\leqslant M((\min F_j)_{j=1}^k)\in\mathcal{S}_{\delta_{I(n)}+1}$. Therefore $(F_j)_{j=1}^k$ is the desired splitting of $E$.

$(v)$ We work by induction on $\delta$ for $\gamma$ held fixed. For $\delta=0$, the result is trivial, since $F$ must be a singleton in this case. We can take $N=\nn$.  

Assume that for $\gamma, \delta<\omega_1$, $N\in[\nn]$ satisfies the conclusions.  Fix $F\in \mathcal{S}_{\delta+1}$ and successive sets $(E_i)_{i\in F}$ such that $\min E_i\geqslant N(i)$ for all $i\in F$.  Write $F=\cup_{j=1}^k F_j$ for some $k\leqslant F_1<\ldots <F_k$, $\varnothing\neq F_j\in\mathcal{S}_\delta$.   By the inductive hypothesis, \[G_j:=\bigcup_{i\in F_j}E_i\in \mathcal{S}_{\gamma+\delta}\] for all $1\leqslant j\leqslant k$.  Since \[k\leqslant \min F_1 \leqslant N(\min F_1)\leqslant \min E_{\min F_1}=\min G_1,\] \[\bigcup_{i\in F}E_i=\bigcup_{j=1}^k G_j\in \mathcal{S}_{\gamma+\delta+1}.\]  

Let $\delta$ be a limit ordinal and assume the result holds for all ordinals less than $\delta$.   Note that $\gamma+\delta$ is a limit ordinal.  Fix $\delta_n\uparrow\delta$ and $\mu_n\uparrow \gamma+\delta$ such that \[\mathcal{S}_\delta=\bigcup_{n=1}^\infty \mathcal{S}_{\delta_n+1}\cap [[n,\infty)]^{<\omega}\] and \[\mathcal{S}_{\gamma+\delta}=\bigcup_{n=1}^\infty \mathcal{S}_{\mu_n+1}\cap [[n,\infty)]^{<\omega}.\]   Fix $I\in[\nn]$ such that for all $n\in\nn$,  $\gamma+\delta_n+1\leqslant \mu_{I(n)}+1$ and \[\mathcal{S}_{\gamma+\delta_n+1}\cap [[I(n),\infty)]^{<\omega}\subset \mathcal{S}_{\mu_{I(n)}+1}.\]  For each $n\in\nn$, fix $N_n\in[\nn]$ such that the conclusion holds for $\delta_n+1$.    Fix $N\in[\nn]$ such that \[N(i)\geqslant \max\{N_1(i), \ldots, N_i(i), I(i)\}.\]    Fix $\varnothing\neq F\in \mathcal{S}_\delta$ and let $n=\min F$. Fix successive sets $(E_i)_{i\in F}\subset \mathcal{S}_\gamma$ such that $\min E_i\geqslant N(i)\geqslant N_n(i)$ for all $i\in F$.  Then $F\in \mathcal{S}_{\delta_n+1}$, so \begin{align*} E:=\bigcup_{i\in F}E_i\in \mathcal{S}_{\gamma+\delta_n+1}\cap [[I(n),\infty)]^{<\omega} \subset \mathcal{S}_{\mu_{I(n)+1}}\cap [[I(n),\infty)]^{<\omega}\subset \mathcal{S}_{\gamma+\delta}. \end{align*}

\end{proof}

\section{Schreier spaces - definitions and basic facts}
\label{sectionschreierspaces}
We recall in this section definition of the Schreier spaces $(X_\xi)_{\xi<\omega_1}$, introduced in \cite{AA}, and their basic properties. We introduce also the definition of a $\varrho$-Schreier pair $\varrho<\omega_1$, central for the main results.

We let  $c_{00}$ denote the space of all finitely-nonzero scalar sequences. We let $(e_i)_{i=1}^\infty$ denote the canonical basis of $c_{00}$. We let $(f_i)_{i=1}^\infty$ denote the coordinate functionals to $(e_i)_{i=1}^\infty$, so that \[f_i(e_j)=\delta_{i,j}=\left\{\begin{array}{ll} 1 & : i = j \\ 0 & : i\neq j.\end{array}\right.\]    For $E\in 2^\nn$, we let $E$ denote the projection on $c_{00}$ given by \[Ex=\sum_{i\in E}f_i(x)e_i.\]    We define the $\xi^{th}$ \emph{Schreier space} $X_\xi$ to be the completion of $c_{00}$ with respect to the norm \[\|x\|_\xi=\sup_{E\in \mathcal{S}_\xi} \|Ex\|_{\ell_1},\ \  x\in c_{00}.\] We note that $X_\xi$ has a lattice structure. When there is need to specify, we let $(e^\xi_i)_{i=1}^\infty$ denote the canonical $c_{00}$ basis considered with the norm $\|\cdot\|_\xi$.  Similarly, we let $(f^\xi_i)_{i=1}^\infty$ denote the coordinate functionals in $X^*_\xi$.   

For $C>0$ and sequences $(x_i)_{i=1}^\infty$, $(y_i)_{i=1}^\infty$ in possibly different Banach spaces, we write $(x_i)_{i=1}^\infty\lesssim_C (y_i)_{i=1}^\infty$ if  $(y_i)_{i=1}^\infty$ $C$-dominates $(x_i)_{i=1}^\infty$. That is, if for all $(a_i)_{i=1}^\infty\in c_{00}$, \[\Bigl\|\sum_{i=1}^\infty a_ix_i\Bigr\|\leqslant C\Bigl\|\sum_{i=1}^\infty a_iy_i\Bigr\|.\]    We write $(x_i)_{i=1}^\infty\lesssim (y_i)_{i=1}^\infty$ if $(x_i)_{i=1}^\infty\lesssim_C (y_i)_{i=1}^\infty$ for some $C$.   We write $(x_i)_{i=1}^\infty \approx (y_i)_{i=1}^\infty$ if  $(x_i)_{i=1}^\infty$ and $(y_i)_{i=1}^\infty$ are equivalent, i.e. $(x_i)_{i=1}^\infty\lesssim (y_i)_{i=1}^\infty$ and $(y_i)_{i=1}^\infty\lesssim (x_i)_{i=1}^\infty$.   

We say a basic sequence $(x_i)_{i=1}^\infty$ is $1$-\emph{right dominant} (resp. $1$-\emph{left dominant}) if whenever $K,L\in[\nn]$ are such that $K(i)\leqslant L(i)$ for all $i\in\nn$ (resp. $K(i)\geqslant L(i)$ for all $i\in\nn$), then $(x_i)_{i\in K}\lesssim_1 (x_i)_{i\in L}$.   

For a Banach space $X$, $x\in X$, and $x^*\in X$, we let $x\otimes x^*$ denote the operator from $X$ into $X$ given by $(x\otimes x^*)(y)=x^*(y)x$.   Similarly, $x^*\otimes x$ denotes the operator from $X^*$ into $X^*$ given by $(x^*\otimes x)(y^*)=y^*(x) x^*$.  

For $C>0$, an ordinal $\xi<\omega_1$, and  a seminormalized sequence $(x_i)_{i=1}^\infty$ in some Banach space $X$, we say that $(x_i)_{i=1}$ is a $C$-$\ell_1^\xi$ \emph{spreading model} if \[\inf\Bigl\{\Bigl\|\sum_{i\in E}a_ix_i\Bigr\|:E\in \mathcal{S}_\xi, \sum_{i\in E}|a_i|=1\Bigr\}\geqslant C.\]  Note that this inequality is trivial when $\xi=0$.    We say a seminormalized sequence $(x_i)_{i=1}^\infty$ is a $C$-$c_0^\xi$ spreading model if \[\sup\Bigl\{\Bigl\|\sum_{i\in E}a_ix_i\Bigr\|:E\in \mathcal{S}_\xi, \max_{i\in E}|a_i|\leqslant 1\Bigr\}\leqslant C.\]   Again, this inequality is trivial if $\xi=0$.   We say that $(x_i)_{i=1}^\infty$ is an $\ell_1^\xi$ spreading model (resp. $c_0^\xi$ spreading model) if it is a $C$-$\ell_1^\xi$ (resp. $C$-$c_0^\xi$) spreading model for some $C>0$.  

 We state below the classical properties of the canonical bases of Schreier spaces and their duals, cf. \cite{AA, BL, GL}.

\begin{proposition} Fix $\xi<\omega_1$. 
\begin{enumerate}[(i)]\item The basis $(e^\xi_i)_{i=1}^\infty$ is normalized, $1$-unconditional, shrinking, and $1$-right dominant.  It is also a 1-$\ell_1^\xi$ spreading model and has no subsequence which is an $\ell_1^{\xi+1}$ spreading model.

\item  The basis $(f^\xi_i)_{i=1}^\infty$ is normalized, $1$-unconditional, 
and $1$-left dominant.  It is also a 1-$c_0^\xi$ spreading model  and has no subsequence which is a $c_0^{\xi+1}$  spreading model. 

\item For $L,M\in[\nn]$ such that $L(i)\leqslant M(i)<L(i+1)$ for all $i\in\nn$, $(e_i^\xi)_{i\in M}\lesssim_2 (e^\xi_i)_{i\in L}$ and $(f^\xi_i)_{i\in L}\lesssim_2 (f^\xi_i)_{i\in M}$. 

\item For $C>0$, any Banach space $X$, and any $1$-unconditional basic sequence $(x_i)_{i=1}^\infty\subset X$, $(x_i)_{i=1}^\infty$ is a $C$-$\ell_1^\xi$ spreading model if and only if $(e^\xi_i)_{i=1}^\infty\lesssim_{C^{-1}}(x_i)_{i=1}^\infty$.  

\item For $C>0$, any Banach space $X$, and any seminormalized sequence $(x_i)_{i=1}^\infty\subset X$, $(x_i)_{i=1}^\infty$ is a $C$-$c_0^\xi$ spreading model if and only if $(x_i)_{i=1}^\infty\lesssim_C (f_i^\xi)_{i=1}^\infty$.   

\end{enumerate}
\label{schreierspacefacts}
\end{proposition}

\begin{proof}

$(i)$ The basis $(e_i^\xi)_{i=1}^\infty$ is a normalized 1-$\ell_1^\xi$ spreading model by definition, 1-unconditional as $\mathcal{S}_\xi$ is hereditary and 1-right dominant as
$\mathcal{S}_\xi$ is spreading. The proof that the basis is shrinking is given in \cite{AA}, short reasoning for $\xi=1$ that could be extended directly to any $\xi<\omega_1$ is presented in \cite[Prop. 3.10]{BL}. The fact that the basis does not admit a subsequence that is an $\ell_1^{\xi+1}$ spreading model follows from Corollary \ref{smallbetanorm} proved in the sequel. 

$(ii)$ The first statement follows from corresponding properties of the basis of $X_\xi$. 
The basis is a 1-$c_0^\xi$ spreading model by unconditionality and the fact that for any $E\in\mathcal{S}_\xi$ the basis $(e_i^\xi)_{i\in E}$ is 1-equivalent to the canonical basis of $\ell_1^{|E|}$. The basis has no subsequence which is a $c_0^{\xi+1}$ spreading model by $(i)$ and duality argument. 

$(iii)$  The statement for subsequences of the basis $(e_i^\xi)_{i=1}^\infty$ follows by the fact that for any set $E=(E(M(i)))_{i\in A}\in\mathcal{S}_\xi\cap [M]^{<\omega}$ we have $(E(L(i)))_{i\in A\setminus \min A}\in\mathcal{S}_\xi$ by the assumption on $M,L$ and the spreading property of $\mathcal{S}_\xi$, thus also $(E(L(i)))_{i\in A}\in \mathcal{A}_2[\mathcal{S}_\xi]$. The statement for  subsequences of the basis $(f_i^\xi)_{i=1}^\infty$ follows by duality. 

$(iv)$ and $(v)$ The properties follows directly by definition. 

\end{proof}







We extend here the notion of an unrestricted Schreier space introduced in \cite[Section 3.2]{BL}. 

Given $\xi<\omega_1$ define $Z_\xi=\{(a_i)_{i=1}^\infty\in\mathbb{R}^\nn: \|(a_i)_{i=1}^\infty\|_\xi<\infty\}$. Obviously $X_\xi$ is a closed subspace of $Z_\xi$, spanned by $(e^\xi_i)_{i=1}^\infty$. Note that the coordinate functionals $(f^\xi_i)_{i=1}^\infty$ extended to $Z_\xi$ have also norm 1. 

Repeating the reasoning in \cite[Section 3.14]{BL} for $p=1$ and  arbitrary $\xi<\omega_1$ (as the standard basis $(e^\xi_i)_{i=1}^\infty$ of $X_\xi$ is bimonotone and shrinking by Prop. \ref{schreierspacefacts} (i)) we obtain  that $Z_\xi$ can be identified with the bidual space $X^{**}_\xi$ in the following way.
\begin{lemma} There is an isometric isomorphism $\Psi: X_\xi^{**}\to Z_\xi$, such that $\Psi\circ \kappa$, where $\kappa$ is the canonical embedding of $X_\xi$ into its bidual space, is the inclusion map $X_\xi\hookrightarrow Z_\xi$.
\label{Zxi}
\end{lemma}

For any $M\in[\nn]$ let $Z_\xi(M)=\{(a_i)_{i=1}^\infty\in Z_\xi: a_i=0 \text{ for any }i\in\nn\setminus M\}$ and $X_\xi(M):=Z_\xi(M)\cap X_\xi=\overline{\mathrm{span}}\{e^\xi_i: i\in M\}$. Obviously $X_\xi(M)$ is 1-complemented in $X_\xi$ by the coordinate projection $P_M$, and $Z_\xi(M)$ is 1-complemented in $Z_\xi$ by $P^{**}_M$. 

\ 

We introduce now the definition of Schreier pairs we shall refer to in the main results.  
\begin{definition} Fix $\xi<\omega_1$ and $\varrho<\omega_1$. 
Let $(x_i)_{i=1}\subset X_\xi$, $(x_i^*)_{i=1}^\infty\subset X_\xi^*$ be normalized block sequences. We say that $((x_i)_{i=1}^\infty, (x_i^*)_{i=1}^\infty)$ is a $\varrho$-Schreier pair if for some $J\in[\nn]$ we have the following.
    \begin{enumerate}[$(i)$]\item $x^*_i(x_j)=0$ for all $i,j\in\nn$ with $i\neq j$, 

\item $\text{Re\ }x^*_i(x_i)=|x^*_i|(|x_i|)$ for all $i\in\nn$, 

\item $\inf_i \text{Re\ } x^*_i(x_i)>0$, 

\item $(x_i)_{i=1}^\infty\approx (e^{\varrho}_i)_{i\in J}$, 

\item $(x^*_i)_{i=1}^\infty\approx (f^{\varrho}_i)_{i\in J}$, 

\item $\sup_n \|\sum_{i=1}^n x_i\otimes x^*_i\|<\infty$. 
\end{enumerate}

\end{definition}

\begin{rem} Let $((x_i)_{i=1}^\infty, (x_i^*)_{i=1}^\infty)\subset X_\xi\times X_\xi^*$ be a $\varrho$-Schreier pair. Then $(x_i)_{i=1}^\infty$ and $(x^*_i)_{i=1}^\infty$ span complemented subspaces of $X_\xi$ and $X_\xi^*$, respectively.

Indeed, it follows from $1$-unconditionality together with the fact that $\inf_i x^*_i(x_i)>0$ that the strong operator topology limit of $\sum_{i=1}^n \frac{x_i\otimes x^*_i}{x^*_i(x_i)}$ defines a bounded projection on the closed span of $(x_i)_{i=1}^\infty$ in $X_\xi$ the adjoint of which is a bounded projection onto the closed span of $(x^*_i)_{i=1}^\infty$ in $X^*_\xi$. 

\label{complemented}
    
\end{rem}

\section{Repeated averages hierarchy}

This section is devoted to the hierarchy of repeated averages, introduced in \cite{AMT}. We recall the definition and the basic properties, and prove Lemma \ref{weaksumming},  which yields estimate of the norm $\|\cdot\|_\xi$ on the sums of vectors in $X_\xi$ generated by repeated averages and in consequence forms a tool for building sequences in $X_\xi$ equivalent to subsequences of the basis of arbitrary $X_\varrho$, $\varrho\in R(\xi)$. 

We define now the repeated averages hierarchy.  For each $\xi<\omega_1$, $M\in[\nn]$, and $n\in\nn$, we will define a function $\mathbb{S}^\xi_{M,n}:\nn\to [0,1]$.   

We define \[\mathbb{S}^0_{M,n}(i)= \left\{\begin{array}{ll} 1 & : i=M(n) \\ 0 & i\in \nn\setminus M(n). \end{array}\right.\] 

Assume that $\mathbb{S}^\xi_{M,n}$ has been defined for each $M\in[\nn]$ and $n\in\nn$.   Fix $M\in[\nn]$.   We define $\mathbb{S}^{\xi+1}_{M,n}$ by induction on $n$.  Let \[M_n=M\setminus\cup_{i=1}^{n-1}\text{supp}(\mathbb{S}^{\xi+1}_{M,i})\] and let \[\mathbb{S}^{\xi+1}_{M,n}= \frac{1}{\min M_n}\sum_{j=1}^{\min M_n} \mathbb{S}^\xi_{M_n,j}.\]    

Assume that $\xi$ is a limit ordinal and $\mathbb{S}^\zeta_{M,n}$ has been defined for all $\zeta<\xi$, $M\in[\nn]$, and $n\in\nn$. Fix $\xi_n\uparrow \xi$ such that \[\mathcal{S}_\xi=\bigcup_{n=1}^\infty \mathcal{S}_{\xi_n+1}\cap [[n,\infty)]^{<\omega}.\]     For any $M\in[\nn]$. We define $\mathbb{S}^\xi_{M,n}$ by induction on $n$. We let \[M_n=M\setminus \bigcup_{i=1}^{n-1}\text{supp}(\mathbb{S}^\xi_{M,i})\] and define \[\mathbb{S}^\xi_{M,n} = \mathbb{S}^{\xi_{\min M_n+1}}_{M_n,1}.\] 

 We group below standard properties of the repeated averages hierarchy.

\begin{lemma} Fix $\xi<\omega_1$.

\begin{enumerate}[(i)]

\item For any $M,N\in[\nn]$ and $m,n\in\nn$ such that $\text{\emph{\supp}}(\mathbb{S}^\xi_{M,m})=\text{\emph{supp}}(\mathbb{S}^\xi_{N,n})$, $\mathbb{S}^\xi_{M,m}=\mathbb{S}^\xi_{N,n}$.  

\item For any $M\in[\nn]$, there exists a unique partition $M=\cup_{i=1}^\infty E_i$ with $E_1<E_2<\ldots$ and $E_n\in MAX(\mathcal{S}_\xi)$ for all $n\in\nn$, and $\text{\emph{supp}}(\mathbb{S}^\xi_{M,n})=E_n$ for all $n\in\nn$. 

\item If $M\in[\nn]$ and $E_1<E_2<\ldots$ is the unique partition of $M$ into successive, maximal members of $\mathcal{S}_\xi$, then for any $n\in\nn$ and $E_n<N\in[\nn]$, $\mathbb{S}^\xi_{E_n\smallfrown  N,1}=\mathbb{S}^\xi_{M,n}$.   

\item For any $M\in[\nn]$, if $E_1<E_2<\ldots$ is the unique partition of $M$ into successive, maximal members of $\mathcal{S}_{\xi+1}$, then there exist $0=p_0<p_1<\ldots$ such that $p_n-p_{n-1}=\min E_n$ and \[\mathbb{S}^{\xi+1}_{M,n}= \frac{1}{\min E_n}\sum_{m=p_{n-1}+1}^{p_n} \mathbb{S}^\xi_{M,m}.\]   

\item For any $M\in [\nn]$, if $\xi$ is a limit ordinal and $\xi_n\uparrow\xi$ is such that \[\mathcal{S}_\xi=\bigcup_{n=1}^\infty \mathcal{S}_{\xi_n+1}\cap [[n,\infty)]^{<\omega}.\] and if $E_1<E_2<\ldots$ is the unique partition of $M$ into successive, maximal members of $\mathcal{S}_\xi$, then for all $n\in\nn$, \[\mathbb{S}^\xi_{M,n}= \mathbb{S}^{\xi_{\min E_n}+1}_{M_n,1}=\frac{1}{\min E_n}\sum_{m=1}^{\min E_n}\mathbb{S}^{\xi_{\min E_n}}_{M_n,m}.\]    

\end{enumerate}
\label{repeatedaveragesfacts}
\end{lemma}

\begin{proof} Most properties or their reformulations are stated  in  \cite{AMT,GL}, we include short justification below.  

First we note that for any $M\in [\nn]$ and $\xi<\omega_1$ we have the following.
\begin{enumerate}[(a)]
\item There is a unique partition $M=\bigcup_{i=1}^\infty E_i$ into successive maximal members of $\mathcal{S}_\xi$ (cf.  \cite[Lemma 2.7]{AG}). 

\item If $E_1<E_2<\dots$ is the unique partition of $M$ into successive maximal members of $\mathcal{S}_{\xi+1}$, then there is a sequence $0=p_0<p_1<\dots$ with $p_n-p_{n-1}=\min E_n$ so that $E_n=\bigcup_{m=p_{n-1}+1}^{p_n}F_m$ for all $n\in\nn$, with $F_1<F_2<\dots$ and $F_m\in MAX (\mathcal{S}_\xi)$ for each $m\in\nn$. 

\item If $\xi$ is limit with $\xi_n\uparrow\xi$ chosen as in the definition of $\mathcal{S}_\xi$ and 
$E_1<E_2<\dots$ is the unique partition of $M$ into successive maximal members of $\mathcal{S}_\xi$, then $E_n\in MAX(\mathcal{S}_{\min E_n+1})$ for any $n\in\nn$. 
\end{enumerate}
It follows straightforward by induction on $\xi<\omega_1$,  that $\text{supp}(\mathbb{S}_{M,n}^\xi)\in MAX (\mathcal{S}_\xi)$ for each $n\in\nn$ and $M=\bigcup_{n=1}^\infty\text{supp}(\mathbb{S}_{M,n}^\xi)$, which together with (a) prove $(i)$. 

The property $(ii)$ is proved easily by induction on $\xi<\omega_1$ using (b) and (c), the property $(iii)$ follows by $(i)$ and $(ii)$. 
The property $(iv)$ follows by $(ii)$, $(b)$ and the definition of $(\mathbb{S}_{M,n}^\xi)_n$, the property $(v)$ - again by the definition of $(\mathbb{S}_{M,n}^\xi)_n$, $(ii)$, (c) and $(iv)$.   

\end{proof}

The next result provides a uniform bound on the norm $\|\cdot\|_\xi$, $\xi<\omega_1$, of vectors $\sum\limits_{n=1}^\infty\sum\limits_{i=1}^\infty\mathbb{S}_{M,n}^\xi(i)e_i^\xi\in\ell_\infty$, for any $M\in[\nn]$. The case of $\xi<\omega$ was treated in \cite[Lemma 2.3]{GL} with the estimate depending on $\xi$.

\begin{lemma} For $\xi<\omega_1$, \[\sup_{M\in[\nn]} \sup\Bigl\{\Bigl\|F\sum_{n=1}^\infty \sum_{i=1}^\infty \mathbb{S}^\xi_{M,n}(i)e_i^\xi\Bigr\|_{\ell_1}:\varnothing\neq F\in \mathcal{S}_\xi\cap [M]^{<\omega}\Bigr\}\leqslant 6.\] 
\label{weaksumming}
\end{lemma}

\begin{proof} We work by induction on $\xi$. The $\xi=0$ case is trivial.  

Assume the result holds for some $\xi$.  Fix $M\in[\nn]$ and $\varnothing\neq F\in \mathcal{S}_{\xi+1}\cap [M]^{<\omega}$.  Write $F=\cup_{j=1}^m F_j$ for some $F_1<\ldots <F_m$, $\varnothing\neq F_j\in \mathcal{S}_\xi$ with $m\leqslant \min F$.  Fix $E_1<E_2<\ldots$ such that $E_i\in MAX(\mathcal{S}_{\xi+1})$ and $\cup_{i=1}^\infty E_i=M$.  For each $n\in\nn$, define \[x_n=\sum_{i=1}^\infty \mathbb{S}^{\xi+1}_{M,n}(i)e_i^\xi.\]  By Lemma \ref{repeatedaveragesfacts}, $\supp(x_n)=E_n$. Moreover, there exist $0=p_0<p_1<\ldots$ such that for all $n\in\nn$, $p_n-p_{n-1}=\min E_n$ and \[x_n=\frac{1}{\min E_n}\sum_{m=p_{n-1}+1}^{p_n} \sum_{i=1}^\infty \mathbb{S}^\xi_{M,m}(i)e_i^\xi.\]  By the inductive hypothesis, $\|F_jx\|_{\ell_1}\leqslant 6/\min E_n$ for all $n\in\nn$ and $1\leqslant j\leqslant m$.   Let \[i_0=\min \{i\in\nn:Fx_i\neq 0\}\] and note that \[m\leqslant \min F\leqslant \max E_{i_0}.\] Therefore for any $i\in \nn$,  \[\|Fx_i\|_{\ell_1} \leqslant \min \Bigl\{\|x_i\|_{\ell_1}, \sum_{j=1}^m \|F_jx_i\|_{\ell_1}\Bigr\} \leqslant \min \{1, 6m/\min E_i\} \leqslant \min \{1, 6\max E_{i_0}/\min E_i\}.\]  As noted in Proposition \ref{schreierfamilyfacts}, for $i\in\nn$, we also have \[\|Fx_{i_0+i}\|_{\ell_1}\leqslant 6\max E_{i_0}/\min E_{i_0+i} \leqslant 6\cdot 2^{1-i}.\]  Since $Fx_i=0$ for all $i<i_0$, it follows that  \begin{align*} \Bigl\|F\sum_{i=1}^\infty x_i\Bigr\|_{\ell_1} & \leqslant \sum_{i=i_0}^{i_0+2}1 + 6\sum_{i=3}^\infty \frac{\max E_i}{\min E_{i_0+i}}=3+6\sum_{i=3}^\infty 2^{1-i}=6.\end{align*} 
Assume $\xi$ is a limit ordinal and the result holds for all $\zeta<\xi$.  Fix $\xi_n\uparrow \xi$ such that \[\mathcal{S}_\xi = \bigcup_{n=1}^\infty \mathcal{S}_{\xi_n+1}\cap [[n,\infty)]^{<\omega}.\] Fix $M\in [\nn]$ and $\varnothing\neq F\in \mathcal{S}_\xi\cap [M]^{<\omega}$.   For each $n\in\nn$, define \[x_n=\sum_{i=1}^\infty \mathbb{S}^\xi_{M,n}(i)e_i^\xi.\]   Let $E_1<E_2<\ldots$ be maximal members of $\mathcal{S}_\xi$ such that $M=\cup_{i=1}^\infty E_i$ and note that $\supp(x_n)=E_n$ for all $n\in\nn$ by Lemma \ref{repeatedaveragesfacts}.  Furthermore, with $M_n=M\setminus \cup_{i=1}^{n-1}E_i$, we have 
\[x_n =\sum_{i=1}^\infty\mathbb{S}^{\xi_{\min E_n}+1}_{M_n,1}(i)e_i^\xi = \frac{1}{\min E_n}\sum_{m=1}^{\min E_n} \sum_{i=1}^\infty\mathbb{S}^{\xi_{\min E_n}}_{M_n,m}(i)e_i^\xi.\] By the inductive hypothesis, $\|x_n\|_{\xi_n}\leqslant 6/\min E_n$ for all $n\in\nn$.   Let \[i_0=\min \{i\in\nn:Fx_i\neq \varnothing\}\] and note that $F\in \mathcal{S}_{\xi_{\min F}+1}\subset \mathcal{S}_{\xi_i}$ for all $i>i_0$.   Therefore for all $i>i_0$,  $\|Fx_i\|_{\ell_1}\leqslant 6/\min E_i\leqslant 6\max E_{i_0}/\min E_i$. Moreover, $\|Fx_i\|_{\ell_1}\leqslant 1$ for all $i\in\nn$ and $Fx_i=0$ for all $i<i_0$.   We can now repeat all calculations from the end of the successor case.   

\end{proof}

The following is an immediate consequence of Lemma \ref{weaksumming}. 

\begin{corollary} For any $\xi<\omega_1$ and $M\in [\nn]$, \[\Bigl(\sum_{i=1}^\infty \mathbb{S}^\xi_{M,n}(i)e_i^\xi\Bigr)_{n=1}^\infty \in 6 B_{\ell_1^w(X_\xi)}.\] 

\label{weaksummingcorollary}
\end{corollary}


\begin{corollary} For any $\beta<\xi$ and $\ee\in (0,1)$, there exists $n=n(\beta,\xi,\ee)\in\nn$ such that for any $M\in[\nn]$ with $n\leqslant M$, \[\Bigl\|\sum_{i=1}^\infty \mathbb{S}^\xi_{M,1}(i)e_i^\beta\Bigr\|_\beta\leqslant \ee.\] 
\label{tail}
\end{corollary}

\begin{proof} If $\xi$ is a successor ordinal, say $\xi=\zeta+1$, then $\zeta\geqslant \beta$.  There exists $m\in\nn$ such that $\mathcal{S}_\beta\cap [[m,\infty)]^{<\omega}\subset \mathcal{S}_\zeta$.   Fix $n\geqslant m$ such that $n>6/\ee$.   For any $M\in[\nn]$ with $n\leqslant M$, by Lemma \ref{repeatedaveragesfacts} and Lemma \ref{weaksumming}, \[\Bigl\|\sum_{i=1}^\infty \mathbb{S}^{\xi}_{M,1}(i)e_i^\beta\Bigr\|_\beta\leqslant \Bigl\|\sum_{i=1}^\infty \mathbb{S}^\xi_{M,1}(i)e_i^\zeta\Bigr\|_\zeta = \frac{1}{\min M}\Bigl\|\sum_{k=1}^{\min M}\sum_{i=1}^\infty \mathbb{S}^\zeta_{M,k}(i)e_i^\zeta\Bigr\|_\zeta\leqslant 6/\min M\leqslant 6/n<\ee.\] 

Suppose that $\xi$ is a limit ordinal and fix $\xi_n\uparrow \xi$ such that $\mathcal{S}_\xi=\bigcup_{i=1}^\infty \mathcal{S}_{\xi_i+1}\cap [[i,\infty)]^{<\omega}$.  Fix $l\in\nn$ such that $\beta<\xi_l$ and $m\in\nn$ such that $\mathcal{S}_\beta\cap[[m,\infty)]^{<\omega}\subset \mathcal{S}_{\xi_l}$.   Note that this implies that $\mathcal{S}_\beta\cap [[m,\infty)]^{<\omega}\subset \mathcal{S}_{\xi_i}$ for all $i\geqslant l$, since $\mathcal{S}_{\xi_i}\subset \mathcal{S}_{\xi_i+1}\subset \mathcal{S}_{\xi_{i+1}}$ for all $i\in\nn$.  Fix $n\geqslant \max\{l,m\}$ such that $6/\ee<n$.  Then for any $M\in[\nn]$ with $n\leqslant M$,  by Lemma \ref{repeatedaveragesfacts} and Lemma \ref{weaksumming}, \[\Bigl\|\sum_{i=1}^\infty \mathbb{S}^\xi_{M,1}(i)e_i^\beta\Bigr\|_\beta \leqslant \Bigl\|\sum_{i=1}^\infty \mathbb{S}^\xi_{M,1}(i)e_i^{\xi_{\min M}}\Bigr\|_{\xi_{\min M}} = \frac{1}{\min M}\Bigl\|\sum_{k=1}^{\min M}\sum_{i=1}^\infty \mathbb{S}^{\xi_{\min M}}_{M,k}(i)e_i^{\xi_{\min M}}\Bigr\|_{\xi_{\min M}}\leqslant 6/\min M<\ee.\]

\end{proof}

We note the following consequence of  Corollary \ref{tail}, which restated in terms of so-called special convex combinations serves as the basic tool in the study of mixed Tsirelson spaces defined by Schreier families of higher order and spaces built on their basis, cf. \cite{AT}. 

\begin{corollary} For any $\beta<\xi$, $\ee\in(0,1)$, and  $M\in[\nn]$, there exist $E\in MAX(\mathcal{S}_\xi)\cap [M]^{<\omega}$ and $x\in S_{\ell_1^+}$ such that $\text{\emph{supp}}(x)=E$ and $\|x\|_\beta<\ee$. 
\label{smallbetanorm}
\end{corollary}

\begin{proof} We take $E=\supp(\mathbb{S}^\xi_{M,1})$ and $x=\sum_{i=1}^\infty \mathbb{S}^\xi_{M,1}(i)e_i^\beta$ for any $M\in[\nn]$ such that $\min M$ is sufficiently large. 
\end{proof}

We conclude the section with a result on isometric representation of $c_0$ in Schreier spaces provided by vectors generated by repeated averages of order $\xi$. The general result concerning repeated averages of any order $\iota<\xi$ will be presented in the final section (Theorem \ref{existence}).

\begin{corollary} For any $\xi<\omega_1$ and $L\in[\nn]$, there exists $M\in[L]$ such that $\bigl(\sum_{i=1}^\infty \mathbb{S}^\xi_{M,n}(i)e_i^\xi\bigr)_{n=1}^\infty$ is isometrically equivalent to the canonical $c_0$ basis. 
\label{isometricc0}
\end{corollary}

\begin{proof} If $\xi=0$, we take $M=L$ and $x_i=e^0_{M(i)}$.  Assume $1\leqslant \xi$.  If $\xi$ is a successor, say $\xi=\beta+1$, let $\mathcal{G}_n=\mathcal{A}_n[\mathcal{S}_\beta]$ for all $n\in\nn$.  If $\xi$ is a limit ordinal, fix $\xi_n\uparrow\xi$ such that $\mathcal{S}_\xi=\bigcup_{n=1}^\infty \mathcal{S}_{\xi_n+1}\cap[[n,\infty)]^{<\omega}$ and let $\mathcal{G}_n=\mathcal{S}_{\xi_n+1}$ for all $n\in\nn$.  Define the norm $|\cdot|_n$ on $c_{00}$ by $|x|_n=\max_{E\in \mathcal{G}_n}\|Ex\|_{\ell_1}$. If $\xi$ is a successor, note that $|\cdot|_n\leqslant n\|\cdot\|_\beta$.

Let $L_0=L$.  Assume that for some $i\in\nn$, $L_1, \ldots, L_{i-1}$ and $E_1<\ldots <E_{i-1}$, $E_j\in MAX(\mathcal{S}_\xi)$ have been chosen such that $E_j$ is an initial segment of $L_j$ for all $1\leqslant j<i$.   If $i=1$, we fix $L_1\in [L_0]$ arbitrary. Otherwise  by Corollary \ref{tail} we fix $L_i\in [L_{i-1}]$ such that $E_{i-1}<L_i$ and for each $1\leqslant k<i$, \[\Bigl|\sum_{j=1}^\infty \mathbb{S}^\xi_{L_i,1}(j)e_j^\xi\Bigr|_{\max E_k}< \frac{\eta_k}{2^k}\]  
where 
$\eta_k=\min_{1\leqslant k<i}\min_{j\in \text{supp}(\mathbb{S}^\xi_{M,k})}\mathbb{S}^\xi_{M,k}(j)$. 
Let $E_i$ be the initial segment of $L_i$ such that $E_i\in MAX(\mathcal{S}_\xi)$. This completes the recursive construction.

Let $M=\cup_{i=1}^\infty E_i$ and fix $\varnothing \neq E\in \mathcal{S}_\xi\cap [M]^{<\omega}$. Note that $\mathbb{S}^\xi_{M,i}=\mathbb{S}^\xi_{L_i,1}$ and $\supp(\mathbb{S}^\xi_{M,i})=E_i$.   Let $i_0=\min \{i\in\nn:E\cap E_i\neq \varnothing\}$ and note that $E\in \mathcal{G}_{\max E_{i_0}}$.    Consider two cases.   

Case $1$: $E_{i_0}\subset E$.  Since $E_{i_0}$ is a maximal member of $\mathcal{S}_\xi$, it follows that $E=E_{i_0}$, so \[\Bigl\|E\sum_{i=1}^\infty \sum_{j=1}^\infty \mathbb{S}^\xi_{M,i}(j)e_j^\xi\Bigr\|_{\ell_1} = \Bigl\|E \sum_{j=1}^\infty \mathbb{S}^\xi_{M,i_0}(j)e_j^\xi\Bigr\|_{\ell_1}=1.\]

Case $2$: $E_{i_0}\not\subset E$.   Then $\|E\sum_{j=1}^\infty \mathbb{S}^\xi_{M,i_0}(j)e_j\Bigr\|_{\ell_1}\leqslant 1-\eta_{i_0}$.   Moreover, since $E\in \mathcal{G}_{\max E_{i_0}}$, \begin{align*}\Bigl\|E\sum_{i=i_0+1}^\infty \sum_{j=1}^\infty \mathbb{S}^\xi_{M,i}(j)e_j^\xi\Bigr\|_{\ell_1} & \leqslant  \sum_{i=i_0+1}^\infty \Bigl|\sum_{j=1}^\infty \mathbb{S}^\xi_{L_i,1}(j)e_j^\xi\Bigr|_{\max E_{i_0}} \\ & \leqslant \eta_{i_0}\sum_{i=i_0+1}^\infty \frac{1}{2^i} \leqslant \eta_{i_0}.\end{align*}  Therefore \[\Bigl\|E\sum_{i=1}^\infty \sum_{j=1}^\infty \mathbb{S}^\xi_{M,i}(j)e_j^\xi\Bigr\|_{\ell_1} \leqslant \Bigl\|E\sum_{j=1}^\infty \mathbb{S}^\xi_{M,i_0}(j)e_j^\xi\Bigr\|_{\ell_1}+\Bigl\|E\sum_{i=i_0+1}^\infty \sum_{j=1}^\infty \mathbb{S}^\xi_{M,i}(j)e_j^\xi\Bigr\|_{\ell_1} \leqslant 1-\eta_{i_0}+\eta_{i_0}=1.\] 

Since $\varnothing\neq E\in \mathcal{S}_\xi\cap [M]^{<\omega}$ was arbitrary, we deduce that $\sup_n \|\sum_{i=1}^n \sum_{j=1}^\infty \mathbb{S}^\xi_{M,i}(j)e_j^\xi\Bigr\|_\xi\leqslant 1$.  By $1$-unconditionality, we deduce that $\bigl(\sum_{i=1}^\infty \mathbb{S}^\xi_{M,n}(i)e_i^\xi\bigr)_{n=1}^\infty\lesssim_1 (e^0_i)_{i=1}^\infty$. The reverse inequality is clear, since $\bigl(\sum_{i=1}^\infty \mathbb{S}^\xi_{M,n}(i)e_i^\xi\bigr)_{n=1}^\infty$ is a normalized, $1$-unconditional basic sequence.

\end{proof}

\section{Auxiliary types of sequences of functionals}

In this section we introduce the notion of $\xi$-flat and $\xi$-approachable sequences of functionals, $\xi<\omega_1$. We prove one of the cases of the main results; i.e. that any flat sequence of functionals in $X_\xi^*$, which very roughly means ``functionals with mainly small coefficients'' contains a subsequence equivalent to the canonical basis of $\ell_1$.  We conclude this section with several pieces providing upper and lower estimates on the norms of subsequences of normalized block sequences in $X^*_\xi$ whose supports exhibit specified upper or lower estimates on their complexity.

We  show first that for $\beta<\xi<\omega_1$, vectors in  $S_{X^*_\xi}$ with small $c_0$ norm  can be well-normed by vectors in $B_{X_\xi}$ with small $X_\beta$ norm. 

\begin{proposition} Fix $\beta<\xi<\omega_1$ and $\ee\in (0,1)$. There exists $l=l(\beta,\xi,\ee)\in\nn$ such that for any  \[x^*\in S_{X_\xi^*}\cap \frac{1}{l}B_{\ell_\infty}\cap \text{\emph{span}}\{f_j^\xi:j\geqslant l\},\] there exists \[x\in S_{X_\xi}\cap \ee B_{X_\beta} \cap \text{\emph{span}}\{e_j^\xi:j\in \text{\emph{supp}}(x^*)\}\] such that $\text{\emph{Re\ }}x^*(x)>1-\ee$. 
\label{flatproposition}
\end{proposition}

\begin{proof} Fix $1<l\in\nn$ so large that $\mathcal{S}_{\beta+1} \cap [[l,\infty)]^{<\omega}\subset \mathcal{S}_\xi$ and $l-1>1/\ee$.  Fix \[x^*\in S_{X^*_\xi}\cap  \frac{1}{l}B_{\ell_\infty}\cap \text{span}\{f_j^\xi:j\geqslant l\}\] and \[y\in B_{X_\xi}\cap \text{span}\{e^\xi_j:j\in \supp(x^*)\}\] such that $\text{Re\ }x^*(y)=1$. Recursively select $E_1, \ldots, E_l\in \mathcal{S}_\beta\cap[[l,\infty)]^{<\omega}$ such that \[\|E_iy\|_{\ell_1} = \max\Bigl\{\|Ey\|_{\ell_1}:E\in \mathcal{S}_\beta\cap [[l,\infty)]^{<\omega}, E\subset \bigcup_{j=1}^{i-1}E_j\Bigr\}.\]   By our choice of $l$ and Lemma \ref{mod1}, \[E:=\bigcup_{i=1}^l E_i\in\mathcal{S}_{\beta+1}^M\cap [[l,\infty)]^{<\omega} = \mathcal{S}_{\beta+1}\cap[[l,\infty)]^{<\omega}\subset \mathcal{S}_\xi.\]    We note that \[1\geqslant \|Ey\|_{\ell_1} \geqslant \sum_{i=1}^l \|E_iy\|_{\ell_1} \geqslant l \|E_ly\|_{\ell_1},\] from which it follows that $\|E_ly\|_{\ell_1}\leqslant 1/l$.    Moreover, \[\|y-Ey\|_\beta =\sup_{F\in\mathcal{S}_\beta}\|(F\setminus E) y \|_{\ell_1} \leqslant \|E_ly\|_{\ell_1}\leqslant 1/l.\]   We compute \[|x^*(Ey)|\leqslant \|x^*\|_{c_0}\|Ey\|_{\ell_1}\leqslant \|x^*\|_{c_0} \|y\|_\xi \leqslant 1/l,\] so \[\text{Re\ }x^*(y-Ey)\geqslant 1-1/l.\] Let \[x=\frac{y-Ey}{\|y-Ey\|_\xi}.\]  Since $\text{Re\ }x^*(y-Ey)\geqslant  1-1/l$, $\|y-Ey\|_\xi\geqslant 1-1/l$.  Then  \[\text{Re\ }x^*(x)\geqslant \text{Re\ }x^*(y-Ey)\geqslant 1-1/l\geqslant 1-\ee\] and \[\|x\|_\beta\leqslant \frac{1/l}{1-1/l}<\ee.\]

\end{proof}

Proposition \ref{flatproposition} suggests the analysis of "size" of coefficients of functionals as a useful tool for classifying sequences of functionals on Schreier spaces. We introduce below mappings $T_\varepsilon, R_\varepsilon$ that split a functional with respect to the "size" of its coefficients. Next we shall define specific classes ($\xi$-flat and $\xi$-approachable)  of sequences of functionals in terms of these mappings. The index of a sequence of functionals $(y_i^*)_{i=1}^\infty$ related to these classes will determine the ordinal $\varrho$ of the basis $(e^\varrho_i)_{i=1}^\infty$ with a subsequence equivalent to a subsequence of $(y_i^*)_{i=1}^\infty$ as we shall prove in the last section. 

For $\ee>0$, we define $T_\ee,R_\ee:c_{00}\to c_{00}$ by \[T_\ee x^*=\sum_{i:|x^*(e_i)|>\ee} x^*(e_i)f_i\] and \[R_\ee x^*=\sum_{i:|x^*(e_i)|\leqslant \ee}x^*(e_i)f_i.\]

\begin{definition}
Let $(y^*_i)_{i=1}^\infty\subset c_{00}$ be a block sequence.   For $\xi<\omega_1$, we say $(y^*_i)_{i=1}^\infty$ is $\xi$-\emph{flat} if \[\underset{\ee>0}{\ \inf\ }\underset{i}{\ \lim{\sup} 
\ }\|R_\ee y^*_i\|_{\xi}>0.\] 

\end{definition}

\begin{corollary} Fix $1\leqslant \xi<\omega_1$ and $\iota\leqslant \xi$.   Suppose that $(y^*_i)_{i=1}^\infty\subset X^*_\xi$ is a $\xi$-flat normalized block sequence.  Then  for any $\eta>0$, there exist a subsequence $(x^*_i)_{i=1}^\infty$ of $(y^*_i)_{i=1}^\infty$ and a normalized block sequence $(x_i)_{i=1}^\infty\subset X_\xi$ such that $((x_i)_{i=1}^\infty,(x_i^*)_{i=1}^\infty)$ is a 0-Schreier pair with $(x_i)_{i\in I}\approx_{1+\eta} (e^0_i)_{i\in \nn}$. 

\label{flatcorollary}
\end{corollary}

\begin{proof} Without loss of generality, we can replace $y^*_i$ with $|y^*_i|$. Fix $0<\vartheta<\underset{\ee>0}{\ \inf\ }\underset{i}{\ \lim\sup\ }{\|R_\varepsilon y^*_i\|_\xi}$.  If $\xi$ is a successor, say $\xi=\gamma+1$, let $\mathcal{G}_n=\mathcal{A}_n[\mathcal{S}_\gamma]$ for all $n\in\nn$.   If $\xi$ is a limit ordinal, fix $\xi_n\uparrow \xi$ such that \[\mathcal{S}_\xi=\bigcup_{n=1}^\infty \mathcal{S}_{\xi_n+1}\cap[[n,\infty)]^{<\omega}\]   and let  $\mathcal{G}_n=\mathcal{S}_{\xi_n+1}$ for each $n\in\nn$. 
For each $n\in\nn$, let $|\cdot|_n$ denote the norm on $c_{00}$ given by $|x|_n=\sup_{E\in \mathcal{G}_n}\|Ex\|_{\ell_1}$.   

We note that for any $m,k\in \nn$ and $\ee>0$, there exist $i>k$ and \[x\in S_{X_\xi}\cap \ee B_{\mathcal{G}_m}\cap \text{span}\{e^\xi_j:j\in \supp(y^*_i)\}\] such that $|y^*_i|(|x|)=\text{Re\ }y^*_i(x)>\vartheta$.    Here, $B_{\mathcal{G}_m}$ denotes the set of $x\in c_{00}$ such that $|x|_m\leqslant 1$.  Indeed, fix $m,k\in\nn$ and $\ee>0$.   Fix  $\ee_0\in (0,1)$ such that \[\vartheta < (1-\ee_0)\underset{\ee>0}{\ \inf\ }\underset{i}{\ \lim\sup\ }\|R_\ee y^*_i\|_\xi.\]   If $\xi=\gamma+1$, let $\beta=\gamma$ and otherwise let $\beta=\xi_m+1$.   By Proposition \ref{flatproposition}, there exists $l=l(\beta,\xi,\min \{\ee_0,\ee/m\})$ such that for any \[x^*\in S_{X_\xi^*}\cap \frac{1}{l}B_{\ell_\infty}\cap \text{span}\{f^\xi_j:j\geqslant l\},\] there exists \[x\in S_{X_\xi}\cap \min \{\ee_0,\ee/m\} B_{X_\beta} \cap \text{span}\{e^\xi_j:j\in \supp(x^*)\}\] such that $\text{Re\ }x^*(x)>1-\ee$.  Of course, if $x^*=|x^*|$, we can replace $x$ with $|x|$ and retain all of these properties. Therefore we can conclude not only that $\text{Re\ }x^*(x)>1-\ee$, but also that $|x^*|(|x|)=\text{Re\ }x^*(x)>1-\ee$.   We can select $i\geqslant \max\{j,k\}$ such that $\|R_{\vartheta/l}y^*_i\|_\xi>\vartheta/(1-\ee_0)$.   Applying the properties of $l$ to \[x^*=\frac{R_{\vartheta/l}y^*_i}{\|R_{\vartheta/l}y^*_i\|_\xi}\in S_{X_\xi^*}\cap \frac{1}{l}B_{\ell_\infty}\cap \text{span}\{f^\xi_j:j\geqslant l\},\] we deduce the existence of \[x\in S_{X_\xi}\cap \frac{\ee}{m}B_{X_\beta}\cap \text{span}\{e^\xi_j:j\in \supp(R_{\vartheta/l}y^*_i)\}\] such that $|x|=x$ and \[y^*_i(x) > (1-\ee_0)\|R_{\vartheta/l}y^*_i\|_\xi\geqslant \vartheta.\]  Note that $|x|_m\leqslant m\|x\|_\beta\leqslant \ee$.  This gives the claim.

Using the claim, we recursively select $I(1)<I(2)<\ldots$ and $x_i\in S_{X_\xi}$ such that, with $x^*_i=y^*_{I(i)}$, for all $i\in\nn$, 

\begin{enumerate}[$(i)$]\item for all $ i<j$, $|x_j|_{\max \supp(x_i)}<\eta/2^{j}$, 

\item  $|x_i|=x_i$, 

\item $\supp(x_i)\subset \supp(x^*_i)$, 

\item $x^*_i(x_i)\geqslant \vartheta$. 
\end{enumerate}

Since $(x_i)_{i=1}^\infty$ is a normalized block sequence, $(e^0_i)_{i=1}^\infty\lesssim_1 (x_i)_{i=1}^\infty$.   Fix $F\in \mathcal{S}_\xi$.   If $Fx_i=0$ for all $i\in\nn$, then $\|F\sum_{i=1}^\infty x_i\|_{\ell_1}=0$.  Otherwise let \[i_0=\min \{i\in\nn:Fx_i\neq 0\}\] and note that $F\in \mathcal{G}_{\max \supp(x_{i_0})}$, so that \begin{align*} \Bigl\|F\sum_{i=1}^\infty x_i\Bigr\|_{\ell_1} & \leqslant \|x_{i_0}\|_\xi + \sum_{j=i_0+1}^\infty |x_j|_{\max \supp(x_{i_0})}\leqslant 1+\eta.\end{align*} By $1$-unconditionality, $(x_i)_{i=1}^\infty\lesssim_{1+\eta}(e^0_i)_{i=1}^\infty$.  By duality, $(f^0_i)_{i=1}^\infty \lesssim_{(1+\eta)/\vartheta} (x^*_i)_{i=1}^\infty$. Of course, $(x^*_i)_{i=1}^\infty\lesssim_1 (f^0_i)_{i=1}^\infty$.

We will show that  $\sup_n \|\sum_{i=1}^n x_i\otimes x^*_i\|\leqslant 1+\eta$, for which it is sufficient to prove the inequality $\sup_n \|\sum_{i=1}^n x^*_i\otimes x_i\|\leqslant 1+\eta$ on the adjoints.   Since $B_{X^*_\xi}$ is the closed, convex hull of \[\Bigl\{\sum_{i\in F}\ee_if^\xi_i:F\in\mathcal{S}_\xi, \max_{i\in F}|\ee_i|=1\Bigr\}\] and $\sum_{i=1}^n x_i^*\otimes x_i$ is a positive operator, it is sufficient to show that \[\Bigl\|\bigl(\sum_{i=1}^n x^*_i\otimes x_i\bigr)(1_F)\Bigr\|\leqslant 1+\eta\] for any $E\in \mathcal{S}_\xi$, where $1_F=\sum_{i\in F}f^\xi_i$.    However, for any $F\in \mathcal{S}_\xi$, \[\Bigl\|\bigl(\sum_{i=1}^n x^*_i\otimes x_i\bigr)(1_F)\Bigr\|\leqslant \sum_{i=1}^n |1_F(x_i)|\|x^*_i\|_\xi=\sum_{i=1}^n \|Fx_i\|_{\ell_1}=\Bigl\|F\sum_{i=1}^n x_i\Bigr\|_{\ell_1}\leqslant 1+\eta.\] 

\end{proof}

For the next results, we introduce some terminology. 

\begin{definition} Let $(x^*_i)_{i=1}^\infty\subset c_{00}$ be a block sequence.  We say $(x^*_i)_{i=1}^\infty$ is $\iota$-\emph{approachable} if there exists $\vartheta>0$ such that for any $\beta<\iota$ and $n,j\in\nn$, there exists $i>j$ such that \[\text{supp}(T_\vartheta x^*_i)\in [\nn]^{<\omega}\setminus \mathcal{A}_n[\mathcal{S}_\beta]'.\]  Every block sequence is considered to be vacuously $0$-approachable.

\end{definition}

\begin{lemma} Fix $1\leqslant \xi<\omega_1$.   Suppose that $(x^*_i)_{i=1}^\infty$ is a normalized block sequence in $X^*_\xi$ that is $\iota$-approachable for some $1\leqslant \iota\leqslant \xi$. 
Then there exist $I\in[\nn]$ and a normalized block sequence $(x_i)_{i=1}^\infty$ in $X_\xi$ such that 

\begin{enumerate}[(i)]\item $x^*_i(x_j)=0$ for all $i,j\in\nn$ with $i\neq j$, 

\item $\text{\emph{Re\ }}x^*_i(x_i)=|x^*_i|(|x_i|)$ for all $i\in I$, 

\item $\inf_i \text{Re\ } x_i^*(x_i)>0$,

\item $\lim_{i\in I}\|x_i\|_\beta=0$ for all $\beta<\iota$, 

\item $\text{\emph{supp}}(x_i)\in \mathcal{S}_\iota\cap \mathcal{S}_\xi$ for all $i\in I$.

\end{enumerate}
\label{lower2}
\end{lemma}

\begin{proof} We consider the limit case first. Fix $(\beta_i)_{i=1}^\infty\subset \iota$ such that $\lim_i \beta_i+1=\iota$. Without loss of generality, we can replace $x^*_i$ with $|x^*_i|$.   Fix $N\in[\nn]$ such that for all $i\in\nn$, \[\bigcup_{j=1}^i \mathcal{S}_{\beta_j}\cap [[N(i), \infty)]^{<\omega} \subset \mathcal{S}_{\beta_i},\]  \[N(i)\geqslant 6i,\] and \[\mathcal{S}_{\beta_i+1}\cap [[N(i),\infty)]^{<\omega}\subset \mathcal{S}_\iota\cap \mathcal{S}_\xi.\]   
By hypothesis, there exists $I\in[\nn]$ such that for all $i\in \nn$, $\min \supp(x^*_{I(i)})\geqslant N(i)$ and  \[\supp(T_\vartheta x^*_{I(i)}) \in [\nn]^{<\omega}\setminus \mathcal{A}_{N(i)}[\mathcal{S}_{\beta_i}]'.\] Then we can fix $E_i\subset \supp(T_\vartheta x^*_{I(i)})$ such that $E_i\in MAX(\mathcal{A}_{N(i)}[\mathcal{S}_{\beta_i}])$. Fix $E_i<M_i\in[\nn]$ and let \[x_{I(i)}=\frac{1}{N(i)}\sum_{m=1}^{N(i)}\sum_{j=1}^\infty \mathbb{S}^{\beta_i}_{E_i\smallfrown M_i,m}(j)e_j^\xi\in \text{span}\{e^\xi_j:j\in E_i\}.\]  Since $\min E_i\geqslant \min \supp(x^*_{I(i)})\geqslant N(i)$, \[E_i\in \mathcal{A}_{N(i)}[\mathcal{S}_{\beta_i}]\cap [[N(i),\infty)]^{<\omega}\subset \mathcal{S}_{\beta_i+1}\cap [[N(i),\infty)]^{<\omega}\subset \mathcal{S}_\iota\cap \mathcal{S}_\xi.\]  By properties of the repeated averages hierarchy and our choice of $N$,  using Lemma \ref{weaksumming}, we obtain that for any $1\leqslant j\leqslant i$, $\|x_{I(i)}\|_{\beta_j}\leqslant \|x_{I(i)}\|_{\beta_i}\leqslant 6/N(i)\leqslant 1/i$.   From this it follows that $\lim_{i\in I}\|x_i\|_\beta=0$ for all $\beta<\iota$.    For $i\in \nn\setminus I$, let $x_i=e^\xi_{\min \supp(x^*_i)}$.   

 If $\iota$ is a successor ordinal, say $\iota=\beta+1$ we repeat the reasoning above with each $\beta_i$ replaced by $\beta$. 

\end{proof}

\begin{lemma} Fix $1\leqslant \xi<\omega_1$.  Suppose that $\iota, \varrho, \delta$ are ordinals such that $\iota+\varrho=\xi$ and $\varrho\geqslant \omega^\delta$. Suppose that $(x^*_i)_{i=1}^\infty$ is a normalized block sequence in $X^*_\xi$ which is neither $\xi$-flat  nor $(\iota+\omega^\delta)$-approachable. 
Then for any $\eta>0$ and $I\in \nn$, there exists $J\in [I]$ such that $(x^*_i)_{i\in J}$ is a $(1+\eta)$-$c_0^{\varrho}$ spreading model. 
\label{upper1}
\end{lemma}

\begin{proof} Note that for any $I\in[\nn]$, $(x^*_i)_{i\in I}$ satisfies the same hypotheses as $(x^*_i)_{i=1}^\infty$. So it is sufficient to prove that for any $\eta>0$, there exists $I\in[\nn]$ such that $(x^*_i)_{i\in I}$ is a $(1+\eta)$-$c_0^\varrho$ spreading model.   

 As $(x_i^*)_{i=1}^\infty$ is not $(\iota+\omega^\delta)$-approachable, for any $\ee\in(0,1)$, there exist $\alpha_\ee<\iota+\omega^\delta$ and $n_\ee, j_\ee\in \nn$ such that for all $i>j_\ee$, \[\text{{supp}}(T_\ee x^*_i)\in \mathcal{A}_{n_\ee}[\mathcal{S}_{\alpha_\ee}].\]  If we select $l\in\nn$ such that $\mathcal{S}_\iota \cap [[l,\infty)]^{<\omega}\subset \mathcal{S}_{\max\{\iota, \alpha_\ee\}}$, we can replace $j_\ee$ with $\max \{l,j_\ee\}$ and $\alpha_\ee$ with $\max\{\iota, \alpha_\ee\}$ and assume that $\alpha_\ee\in [\iota, \iota+\omega^\delta)$ for all $\ee\in (0,1)$.  We can then write $\alpha_\ee=\iota+\beta_\ee$.  Therefore for all $\ee\in (0,1)$, there exist $\beta_\ee<\omega^\delta$ and $n_\ee, j_\ee\in\nn$ such that for all $i>j_\ee$, \[\text{supp}(T_\ee x^*_i)\in \mathcal{A}_{n_\ee}[\mathcal{S}_{\iota+\beta_\ee}].\] 
Note that $T_1 x^*_i=0$ for all $i\in\nn$, since $\{i\in\nn:|x^*(e_i)|>1\}=\varnothing$ for any $x^*\in S_{X^*_\xi}$.   For each $\ee\in (0,1)$, let \[\nu(\ee)=\underset{i}{\ \lim\sup\ }\|R_\ee x^*_i\|_\xi.\]   As $(x_i^*)_{i=1}^\infty$ is not $\xi$-flat, we can pick  $(\ee_i)_{i=1}^\infty\subset (0,1)$ such that, with $\nu_i=\nu(\ee_i)$, $\sum_{i=1}^\infty \nu_i+\ee_i<\eta/2$.   Fix a sequence $(\beta_i)_{i=1}^\infty\subset \omega^\delta$ and $N,J\in[\nn]$ such that for all $i\in\nn$ and $j\geqslant J(i)$, $\|R_{\ee_i}x^*_j\|<\nu_i+\ee_i$ and for each $1\leqslant k\leqslant i$, \[\supp(T_{\ee_k}x^*_j)\in \mathcal{A}_{N(k)}[\mathcal{S}_{\iota+\beta_k}].\]    By replacing $(x^*_i)_{i=1}^\infty$ with $(x^*_i)_{i\in J}$, we can assume that $J(i)=i$ for all $i\in\nn$.   

Let $\ee_0=1$.   For each $i,j\in\nn$, let \[E_{i,j}=\left\{\begin{array}{ll} \supp(T_{\ee_i}x_j) & : i\leqslant j \\ \varnothing & : i > j.\end{array}\right. \] Note that for each $i\in\nn$, $E_{i,1}<E_{i,2}<\ldots$, $i\leqslant E_{i,j}\in 
\mathcal{A}_{N(i)}[\mathcal{S}_{\iota+\beta_i}]$.   By Proposition \ref{combinatorics}, there exists $I\in[\nn]$ such that \[\bigcup_{i=1}^\infty \Bigl\{\bigcup_{j\in F}E_{i,I(j)}:F\in \mathcal{S}_\varrho\cap [[i,\infty)]^{<\omega}\Bigr\}\subset \mathcal{S}_{\iota+\varrho}=\mathcal{S}_\xi.\]   We will show that $(x^*_i)_{i\in I}$ is a $(1+\eta)$-$c_0^\varrho$ spreading model. By $1$-unconditionality, it is sufficient to show that for any $\varnothing\neq F\in \mathcal{S}_\varrho$, $\|\sum_{j\in F}x_{I(j)}\|_\xi\leqslant 1+\eta$. To that end, fix $\varnothing\neq F\in \mathcal{S}_\varrho$ and note that \begin{align*}\Bigl\|\sum_{j\in F}x^*_{I(j)}\Bigr\|_\xi & \leqslant \sum_{j=1}^\infty \|R_{\ee_j}x^*_{I(j)}\|_\xi + \sum_{i=1}^\infty \Bigl\|\sum_{j\in F\cap [i,\infty)} (T_{\ee_i}x^*_{I(j)}-T_{\ee_{i-1}}x^*_{I(j)})\Bigr\|_\xi \\ & \leqslant \sum_{j=1}^\infty \nu_j+\ee_j + \sum_{j=1}^\infty \ee_{j-1} < \frac{\eta}{2}+1+\sum_{j=1}^\infty \ee_j<1+\eta.\end{align*} Here we have decomposed $x^*_{I(j)}=R_{\ee_j}x^*_{I(j)}+\sum_{i=1}^j (T_{\ee_i}x^*_{I(j)}-T_{\ee_{i-1}}x^*_{I(j)})$. We have also used the fact that since $I(j)\geqslant j$, $\|R_{\ee_j}x^*_{I(j)}\|\leqslant \nu_j+\ee_j$. Finally, we used the fact that for each $i\in\nn$, \[\Bigl\|\sum_{j\in F\cap [i,\infty)} (T_{\ee_i}x^*_{I(j)}-T_{\ee_{i-1}}x^*_{I(j)})\Bigr\|_{c_0}\leqslant \ee_{i-1}\] and \[\supp\Bigl(\sum_{j\in F}(T_{\ee_i}x^*_{I(j)}-T_{\ee_{i-1}}x^*_{I(j)})\Bigr)\subset \bigcup_{j\in F\cap [i,\infty)}E_{i,I(j)}\in \mathcal{S}_\xi,\] so \[\Bigl\|\sum_{j\in F\cap [i,\infty)} (T_{\ee_i}x^*_{I(j)}-T_{\ee_{i-1}}x^*_{I(j)})\Bigr\|_\xi\leqslant \ee_{i-1}.\]

\end{proof}

\section{Equivalence of block sequences}

This section is devoted to the proof of main theorems of the paper. Theorems \ref{final1} and \ref{final2} provide a detailed version of Theorem \ref{main1}, indicating the order $\varrho\in R(\xi)$ of resulting $\varrho$-Schreier pair. The second part of the section is devoted to the existence of $\varrho$-Schreier pairs, $\varrho\in R(\xi)$, in $X_\xi$ and $X_\xi^*$ provided by Theorem \ref{existence}, which together with Corollaries \ref{moveup1} and \ref{moveup2} prove Theorem \ref{main2}.

 The following result generalizes \cite[Prop. 4.3]{GL}, where the case of $\xi<\omega$ was studied. 

\begin{theorem} Let $(y_i)_{i=1}^\infty\subset X_\xi$ be a normalized block sequence.

\[\iota=\max\{\gamma\in I(\xi):(\forall \beta<\gamma)(\underset{i}{\lim\sup} \|y_i\|_\beta=0)\}\] and let $\varrho=\xi-\iota$.  Then there exist  a subsequence $(x_i)_{i=1}^\infty$ of $(y_i)_{i=1}^\infty$,  and a  normalized block sequence $(x^*_i)_{i=1}^\infty\subset X_\xi^*$ such that $((x_i)_{i=1}^\infty, (x_i^*)_{i=1}^\infty)$ is a $\varrho$-Schreier pair. In particular, $(x_i)_{i=1}^\infty$ and $(x^*_i)_{i=1}^\infty$ span complemented subspaces of $X_\xi$ and $X_\xi^*$, respectively. 

\label{final1}

\end{theorem}

\begin{proof} 

We first prove the following.

\begin{lemma} Fix $1\leqslant \xi<\omega_1$ and $\iota\leqslant \xi$. Let $\varrho=\xi-\iota$.    Suppose that $(x_i)_{i=1}^\infty\subset X_\xi$ and $(x^*_i)_{i=1}^\infty\subset X^*_\xi$ are normalized block sequences such that 

\begin{enumerate}[(i)]\item $x^*_i(x_j)=0$ for all $i,j\in\nn$ with $i\neq j$, 

\item $\text{\emph{Re\ }}x^*_i(x_i)=|x^*_i|(|x_i|)$ for all $i\in\nn$, 

\item $\inf_i \text{\emph{Re}\ }x^*_i(x_i)>0$, 

\item $\lim_i \|x_i\|_\beta=0$ for all $\beta<\iota$, 

\item $\text{\emph{supp}}(x^*_i)\in \mathcal{S}_\iota$ for all $i\in\nn$. 

\end{enumerate}
Then there exist $I\in[\nn]$ such that $((x_i)_{i\in I}, (x_i^*)_{i\in I})$ is a $\varrho$-Schreier pair. 

\label{lastlegs1}
\end{lemma}

\begin{proof}[Proof of Lemma \ref{lastlegs1}] Without loss of generality, we can replace $x_i$ with $|x_i|$ and $x^*_i$ with $|x^*_i|$ and assume $(x_i)_{i=1}^\infty\subset X_\xi$ and $(x^*_i)_{i=1}^\infty \subset X^*_\xi$.  

Let 
$\vartheta=\inf_i x^*_i(x_i)$. By Proposition \ref{combinatorics}, there exists $L\in[\nn]$ such that \[\Bigl\{\bigcup_{i\in F} \supp(x^*_{L(i)}): F\in  \mathcal{S}_\varrho\Bigr\}\subset \mathcal{S}_\xi.\]   From this it follows that $(x^*_i)_{i\in L}$ is a $1$-$c_0^\varrho$ spreading model, since $\|\sum_{i\in F}x^*_{L(i)}\|_{c_0}\leqslant 1$ and $\supp(\sum_{i\in F}x^*_{L(i)})\in \mathcal{S}_\xi$ for all $F\in \mathcal{S}_\varrho$.  By duality, using $(i)$ and $(ii)$, it follows that $(x_i)_{i\in L}$ is a $\vartheta^{-1}$-$\ell_1^\varrho$-spreading model.  By passing to a subsequence and relabeling, we can assume $L=\nn$.  From this and $1$-unconditionality of $(x_i)_{i=1}^\infty$, it follows that $(e^\varrho_i)_{i=1}^\infty \lesssim_{\vartheta^{-1}}(x_i)_{i\in K}$ and $(x^*_i)_{i\in K}\lesssim_1 (f^\varrho_i)_{i=1}^\infty$ for all $K\in[\nn]$ (cf. Proposition \ref{schreierspacefacts}).  

If $\iota$ is a successor, say $\iota=\beta+1$, then we define \[\mathcal{G}_n=\mathcal{A}_n[\mathcal{S}_\beta]\] for all $n\in\nn$. If $\iota$ is a limit ordinal, then we fix $\beta_n\uparrow \iota$ such that $\mathcal{S}_\iota=\bigcup_{n=1}^\infty\mathcal{S}_{\beta_n+1}\cap [[n,\infty)]^{<\omega}$ and define $\mathcal{G}_n=\mathcal{S}_{\beta_n+1}$.   Let $|\cdot|_n$ be the norm on $c_{00}$ given by $|x|_n=\sup_{E\in \mathcal{G}_n}\|Ex\|_{\ell_1}$.  By $(iii)$, $\lim_i |x_i|_n=0$ for all $n\in\nn$.  This is because $|\cdot|_n\leqslant n\|\cdot\|_\beta$ and $\beta<\iota$ in the successor case, and  $|\cdot|_n=\|\cdot\|_{\beta_n+1}$ and $\beta_n+1<\iota$ in the limit case.    Fix $\eta>0$. Fix $K\in[L]$ such that for all $1\leqslant i<j$, $i,j\in\nn$, \[|x_{K(j)}|_{\max\supp(x_{K(i)})} < \frac{\eta}{\max\supp(x_{K(i)})2^{i+j}}.\]    By Proposition \ref{combinatorics}, there exists $M\in[\nn]$ such that for any $E\in \mathcal{S}_\xi$, there exist $E_1<\ldots <E_k$, $\varnothing\neq E_i\in \mathcal{S}_\iota$, such that $M((\min E_i)_{i=1}^k)\in \mathcal{S}_\varrho$.  Let $N(i)=M(\max\supp(x_{K(i)}))$ for all $i\in\nn$.   We claim that $(x_i)_{i\in K}\lesssim_{1+\eta}(e^\varrho_i)_{i\in N}$.  To see this, fix $0 \neq (a_i)_{i=1}^\infty\in c_{00}$ and a minimal 
$E\in \mathcal{S}_\xi$ such that \[\Bigl\|E\sum_{i=1}^\infty a_ix_{K(i)}\Bigr\|_{\ell_1}=\Bigl\|\sum_{i=1}^\infty a_ix_{K(i)}\Bigr\|_\xi.\]  We can fix $E_1<\ldots <E_k$, $\varnothing\neq E_i\in \mathcal{S}_\iota$ such that $E=\cup_{i=1}^k E_i$ and $M((\min E_i)_{i=1}^k)\in \mathcal{S}_\varrho$.    For each $1\leqslant j\leqslant k$, let \[\phi(j)=\min \{i\in\nn:E_jx_{K(i)}\neq 0\}.\] By minimality 
of $E$, each of these sets is non-empty.    Let $A=\{\phi(j):1\leqslant j\leqslant k\}$ and for $i\in A$, let $B_i=\phi^{-1}(\{i\})$ and note that $|B_i|\leqslant |\supp(x_{K(i)})|\leqslant \max\supp(x_{K(i)})$.  Note that $G:=M((\max\supp(x_{K(i)}))_{i\in A})$ is a spread of $M((\min E_{\max B_i})_{i\in A})\in \mathcal{S}_\varrho$, so $G\in \mathcal{S}_\varrho$.  Therefore, since $E_l\in \mathcal{G}_{\max\supp(x_{K(i)})}$ whenever $l\in B_i$,  \begin{align*}\Bigl\|\sum_{i=1}^\infty a_i x_{K(i)}\Bigr\|_\xi & = \Bigl\|E\sum_{i=1}^\infty a_ix_{K(i)}\Bigr\|_{\ell_1} \\ & \leqslant \sum_{i\in A}|a_i|+\|(a_i)_{i=1}^\infty\|_{c_0}\sum_{i\in A}\sum_{l\in B_i} \sum_{j=i+1}^\infty \|E_lx_{I(j)}\|_{\ell_1} \\ & \leqslant \sum_{i\in A}|a_i|+\eta \|(a_i)_{i=1}^\infty\|_{c_0}\sum_{i=1}^\infty \frac{|B_i|}{\max\supp(x_{I(i)})2^{i+j}} \\ & \leqslant \sum_{i\in A}|a_i|+\eta \|(a_i)_{i=1}^\infty\|_{c_0}\sum_{i=1}^\infty \sum_{j=1}^\infty \frac{1}{2^{i+j}} \\ & = \Bigl\|G\sum_{i=1}^\infty a_ie^\varrho_{N(i)}\Bigr\|_{\ell_1}+\eta \|(a_i)_{i=1}^\infty\|_{c_0}\leqslant (1+\eta)\Bigl\|\sum_{i=1}^\infty a_ie_{N(i)}^\varrho\Bigr\|_\varrho. \end{align*}
Fix $J\in[\nn]$ such that $J(i+1)>N(J(i))$ for all $i\in\nn$. Let $I(i)=K(J(i))$ and $P(i)=N(J(i))$.   Then \[(e^\varrho_i)_{i\in J}\lesssim_{\vartheta^{-1}} (x_i)_{i\in I}\lesssim_{1+\eta}(e^\varrho_i)_{i\in P}\lesssim_2 (e^\varrho_i)_{i\in J}\] and \[(f^\varrho_i)_{i\in J}\lesssim_2 (f^\varrho_i)_{i\in P} \lesssim_{\vartheta^{-1}} (x_i^*)_{i\in I}\lesssim_{1+\eta}\lesssim_1 (f^\varrho_i)_{i\in J}.\] We have used duality for the relation $(f^\varrho_i)_{i\in P} \lesssim_{\vartheta^{-1}} (x_i^*)_{i\in I}$.

We conclude by showing that $\sup_n \|\sum_{i=1}^n x_{I(i)}\otimes x^*_{I(i)}\|<2(1+\eta)$. To that end, fix $x\in B_{X_\xi}$. Observe that \begin{align*} \Bigl\|\Bigl(\sum_{i=1}^n x_{I(i)}\otimes x^*_{I(i)}\Bigr)(x)\Bigr\| & \leqslant (1+\eta)\Bigl\|\sum_{i=1}^n |x^*_{I(i)}(x)|e^\varrho_{P(i)}\Bigr\|_\varrho \leqslant 2(1+\eta)\Bigl\|\sum_{i=1}^n |x^*_{I(i)}(x)|e^\varrho_{J(i)}\Bigr\|_\varrho \\ & = 2(1+\eta) \sum_{i\in E}|x^*_{I(i)}(x)|\end{align*} for some $E\in[\nn]^{<\omega}$ such that $F=J(E)\in \mathcal{S}_\varrho$. By our choice of $L$ and since $K\in[L]$, \[\bigcup_{i\in F}\supp(x^*_{K(i)})=\bigcup_{i\in E}\supp(x^*_{I(i)})\in \mathcal{S}_\xi,\] so that 

\[\sum_{i\in F}|x^*_{K(i)}| \in B_{X^*_\xi}.\]  Therefore \[\sum_{i\in E}|x^*_{I(i)}(x)| \leqslant \Bigl(\sum_{i\in F}|x^*_{K(i)}|\Bigr)(|x|) \leqslant 1\cdot 1=1. \]  

\end{proof}

We now return to the proof of Theorem \ref{final1}. We pick  a subsequence $(x_i)_{i=1}^\infty$ and a normalized block sequence $(x_i^*)_{i=1}^\infty\subset X^*_\xi$ satisfying properties $(i)$-$(v)$ of Lemma \ref{lastlegs1}. 
Indeed, without loss of generality, we can replace $x_i$ with $|x_i|$.    Fix  $m\in\nn$ such that $\mathcal{S}_\iota\cap [[m,\infty)]^{<\omega}\subset \mathcal{S}_\xi$.   Fix a subsequence $(x_i)_{i=1}^\infty$ of $(y_i)_{i=1}^\infty$ such that $\min \supp(x_1)\geqslant m$ and $\inf_i \|x_i\|_\iota>0$.    For each $i\in\nn$, fix $E_i\in \mathcal{S}_\iota$ such that $E_i\subset \supp(x_i)$ and $\|E_ix_i\|_{\ell_1}=\|x_i\|_\iota$. Letting $x^*_i=\sum_{j\in E_i}f^\xi_j$ ends the construction. Now we apply Lemma \ref{lastlegs1} and Remark \ref{complemented} to finish the proof. 

\end{proof}

\begin{theorem} Let $(y^*_i)_{i=1}^\infty \subset X_\xi^*$ be a normalized block sequence. If $(y^*_i)_{i=1}^\infty$ is $\xi$-flat, let $\iota=\xi$, and otherwise let \[\iota=\max\{\beta\in I(\xi):(y^*_i)_{i=1}^\infty \text{\ is\ }\beta\text{-approachable}\}.\] Let $\varrho=\xi-\iota$.   Then there exist a subsequence $(x^*_i)_{i=1}^\infty$ of $(y_i^*)_{i=1}^\infty$,  and a normalized block sequence $(x_i)_{i=1}^\infty\subset X_\xi$ such that $((x_i)_{i=1}^\infty, (x_i^*)_{i=1}^\infty)$ is a $\varrho$-Schreier pair. In particular, $(x_i)_{i=1}^\infty$ and $(x^*_i)_{i=1}^\infty$ span complemented subspaces of $X_\xi$ and $X_\xi^*$, respectively. 
\label{final2}
\end{theorem}

\begin{proof} 

We first prove the following. 

\begin{lemma} Fix $1\leqslant \xi<\omega_1$ and $\iota\leqslant \xi$.  Let $\varrho=\xi-\iota$. Suppose that $(x_i)_{i=1}^\infty\subset X_\xi$ and $(x^*_i)_{i=1}^\infty\subset X^*_\xi$ are normalized block sequences such that 

\begin{enumerate}[(i)]
\item $x^*_i(x_j)=0$ for all $i,j\in\nn$ with $i\neq j$, 

\item $\text{\emph{Re\ }}x^*_i(x_i)=|x^*_i|(|x_i|)$ for all $i\in\nn$, 

\item $\inf_i \text{\emph{Re}\ }x^*_i(x_i)>0$, 

\item $\lim_i \|x_i\|_\beta=0$ for all $\beta<\iota$, 

\item $\text{\emph{supp}}(x_i)\in \mathcal{S}_\iota$ for all $i\in\nn$, 

\item $(x^*_i)_{i=1}^\infty$ is a $c_0^{\varrho}$-spreading model. 

\end{enumerate}
Then there exist $I\in[\nn]$ such that $((x_i)_{i\in I}, (x_i^*)_{i\in I})$ is a $\varrho$-Schreier pair. 

\label{lastlegs2}
\end{lemma}

Note that if $\iota=\xi$, then the requirement that $(x^*_i)_{i=1}^\infty$ in the preceding lemma is a $c_0^{\xi-\iota}$ spreading model is trivial, since it simply means that $(x^*_i)_{i=1}^\infty$ is seminormalized.

\begin{proof}[Proof of Lemma \ref{lastlegs2}]  Without loss of generality, we can replace $x_i$ with $|x_i|$ and $x^*_i$ with $|x^*_i|$ and assume $(x_i)_{i=1}^\infty\subset X_\xi$ and $(x^*_i)_{i=1}^\infty \subset X^*_\xi$.  

Let 
$\vartheta=\inf_i x^*_i(x_i)$.    Since $(x^*_i)_{i=1}^\infty$ is a $c_0^\varrho$ spreading model, there exists $C\geqslant 1$ such that for any $K\in[\nn]$, $(x^*_i)_{i\in K}\lesssim_C (f^\varrho_i)_{i=1}^\infty$.  By duality, for any $K\in[\nn]$, $(e^\varrho_i)_{i=1}^\infty\lesssim_{C/\vartheta} (x_i)_{i\in K}$.    

We fix $\eta>0$ and choose $L,K,M,N,I,J\in[\nn]$  as in the proof of Lemma \ref{lastlegs1} and deduce that \[(e^\varrho_i)_{i\in J}\lesssim_{C/\vartheta}(x_i)_{i\in I}\lesssim_{1+\eta}(e^\varrho_i)_{i\in P}\lesssim_2 (e^\varrho_i)_{i\in J} \]  and \[(f^\varrho_i)_{i\in J}\lesssim_2 (f^\varrho_i)_{i\in P} \lesssim_{(1+\eta)/\vartheta}(x^*_i)_{i\in I} \lesssim_C(f^\varrho_i)_{i\in J}.\]

It remains to show that $\sup_n \|\sum_{i=1}^n x_{I(i)}\otimes x^*_{I(i)}\|\leqslant 2C+\eta$.   It is sufficient to prove the inequality $\sup_n \|\sum_{i=1}^n x^*_{I(i)}\otimes x_{I(i)}\|\leqslant 2C+\eta$ on the adjoints. Since $\sum_{i=1}^n x^*_{I(i)}\otimes x_{I(i)}$ is a positive operator and $B_{X^*_\xi}$ is the closed, convex hull of the set $\{\sum_{i\in F}\ee_i f^\xi_i:F\in \mathcal{S}_\xi, |\ee_i|\leqslant 1\}$, it is sufficient to prove that for any $\varnothing \neq E\in \mathcal{S}_\xi$,  \[\Bigl\|\Bigl(\sum_{i=1}^n x^*_{I(i)}\otimes x_{I(i)}\Bigr)(1_E)\Bigr\|\leqslant 2C+\eta,\] where $1_E=\sum_{i\in E}f^\xi_i$.   To that end, fix $\varnothing\neq E\in \mathcal{S}_\xi$.  Write $E=\cup_{l=1}^k E_l$ such that $E_1<\ldots <E_l$, $\varnothing \neq E_l\in \mathcal{S}_\iota$, $M((\min E_l)_{l=1}^k)\in \mathcal{S}_\varrho$.    As in the proof of Lemma \ref{lastlegs1}, we define $\phi$, $A$, and $B_i$, $i\in A$ and estimate \begin{align*}\Bigl\|\Bigl(\sum_{i=1}^n x^*_{I(i)}\otimes x_{I(i)}\Bigr)(1_E)\Bigr\|_\xi & =\Bigl\|\sum_{i\in A} \|Ex_{I(i)}\|_{\ell_1} x^*_{I(i)} + \sum_{i\in A}\sum_{l=1}^k \sum_{j=i+1}^\infty \|E_lx_{I(j)}\|_{\ell_1}x^*_{I(j)}\Bigr\|_\xi \\ & \leqslant \Bigl\|\sum_{i\in A}x^*_{I(i)}\Bigr\|_\xi+\sum_{i\in A} \sum_{l=1}^k\sum_{j=i+1}^\infty \|E_lx_{I(j)}\|_{\ell_1} \\ & \leqslant \Bigl\|\sum_{i\in A}x^*_{I(i)}\Bigr\|_\xi+\eta\sum_{i=1}^\infty\sum_{j=1}^\infty \frac{1}{2^{i+j}} \\ & \leqslant \Bigl\|\sum_{i\in A}x^*_{I(i)}\Bigr\|_\xi +\eta.\end{align*}  Since \[G:=M((\max\supp(x_{I(i)}))_{i\in A})\] is a spread of $M((\min \supp(E_{\max B_i}))_{i\in A})\in \mathcal{S}_\varrho$, it follows that $G\in \mathcal{S}_\varrho$.    Since $G=P(A)$, \[\Bigl\|\sum_{i\in A}x^*_{I(i)}\Bigr\|_\xi\leqslant C\Bigl\|\sum_{i\in A} f^\varrho_{J(i)}\Bigr\|_\varrho\leqslant 2C\Bigl\|\sum_{i\in A}f^\varrho_{P(i)}\Bigr\|_\varrho\leqslant 2C.\]  Therefore \[\Bigl\|\Bigl(\sum_{i=1}^n x^*_{I(i)}\otimes x_{I(i)}\Bigr)(1_E)\Bigr\|\leqslant 2C+\eta.\]

\end{proof}

We now return to the proof of Theorem \ref{final2}. If $(y^*_i)_{i=1}^\infty$ is $\xi$-flat, we apply Corollary \ref{flatcorollary}. Otherwise we  pick a subsequence $(x^*_i)_{i=1}^\infty$ and a normalized block sequence $(x_i)_{i=1}^\infty\subset X_\xi$ satisfying properties $(i)$-$(vi)$ of Lemma \ref{lastlegs2}.  
Indeed,  as $(y_i^*)_{i=1}^\infty$ is $\iota$-approachable, we first apply Lemma \ref{lower2} to get sequences satisfying $(i)$-$(v)$ of Lemma \ref{lastlegs2}.   If $\iota=\xi$, then $(vi)$ is trivially satisfied, otherwise $(y_i^*)_{i=1}^\infty$ is not $(\iota+\omega^\delta)$-approachable for some ordinal $\delta$ with $\iota+\omega^\delta\leqslant \xi$, and we apply  Lemma \ref{upper1} to get a further subsequence satisfying $(vi)$ of Lemma \ref{lastlegs2}. Now we apply Lemma \ref{lastlegs2} and Remark \ref{complemented} to end the proof.

\end{proof}

\begin{corollary} Fix $\xi<\omega_1$. 
Assume that $\varrho
<\omega_1$ and  $(y_i)_{i=1}^\infty$ is a normalized block sequence in $X_\xi$ which is equivalent to a subsequence of the basis of $X_\varrho$. Then $\varrho\in R(\xi)$
and there exist a subsequence $(x_i)_{i=1}^\infty$ of $(y_i)_{i=1}^\infty$ and  a normalized block sequence $(x^*_i)_{i=1}^\infty$ in $X^*_\xi$ such that $((x_i)_{i=1}^\infty, (x_i^*)_{i=1}^\infty)$ is a $\varrho$-Schreier pair. 
\label{moveup1}
\end{corollary}

\begin{proof} By Theorem \ref{final1}, there exist $\gamma\in R(\xi)$, a subsequence $(x_i)_{i=1}^\infty$ of $(y_i)_{i=1}^\infty$, and a normalized block sequence $(x^*_i)_{i=1}^\infty$ in $X^*_\xi$ such that $((x_i)_{i=1}^\infty, (x_i^*)_{i=1}^\infty)$ is a $\gamma$-Schreier pair. Then $(x_i)_{i=1}^\infty$ is equivalent to subsequences of the bases of both $X_\gamma$ and $X_\varrho$, from which it follows that $\gamma=\varrho$. 

\end{proof}

\begin{corollary} Fix $\xi<\omega_1$. 
Assume that $\varrho<\omega_1$ and  $(y_i^*)_{i=1}^\infty$ is a normalized block sequence in $X_\xi^*$ which is equivalent to a subsequence of the basis of $X_\varrho^*$. Then $\varrho\in R(\xi)$ and there exist a normalized block sequence  $(x_i)_{i=1}^\infty\subset X_\xi$,  and a subsequence $(x^*_i)_{i=1}^\infty$ of $(y^*_i)_{i=1}^\infty$ such that $((x_i)_{i=1}^\infty, (x_i^*)_{i=1}^\infty)$ is a $\varrho$-Schreier pair. 
\label{moveup2}
\end{corollary}

\begin{proof} By Theorem \ref{final2}, there exist $\gamma\in R(\xi)$, a normalized block sequence  $(x_i)_{i=1}^\infty\subset X_\xi$, and a subsequence $(x^*_i)_{i=1}^\infty$ of $(y^*_i)_{i=1}^\infty$ such that $((x_i)_{i=1}^\infty, (x_i^*)_{i=1}^\infty)$ is a $\gamma$-Schreier pair. But then $(x^*_i)_{i=1}^\infty$ is equivalent to subsequences of the bases of both $X^*_\gamma$ and $X^*_\varrho$, from which it follows that $\gamma=\varrho$.

\end{proof}

 We proceed now to results on existence of block sequences in $X_\xi$ and $X_\xi^*$ equivalent to subsequences of $(e_i^\varrho)_{i=1}^\infty$ for $\varrho\in R(\xi)$. The main tool is provided by repeated averages hierarchy. 

The next result should be compared with Corollary \ref{isometricc0} concerning the case of $\iota=\xi$, in which case stronger property holds true, namely the constructed sequence is isometrically equivalent to the canonical basis of $c_0$. The case of $\xi<\omega$ was considered in \cite[Lemma 3.10]{GL}.

\begin{theorem}
Fix $\xi<\omega_1$, and $\iota\leqslant \xi$, and let   $\varrho=\xi-\iota$. Then for any $L\in[\nn]$, there exist $M\in[L]$ such that, with $x_i=\sum_{j=1}^\infty \mathbb{S}^\iota_{M,i}(j)e_j^\xi$ and $x^*_i=\sum_{j\in \text{\emph{supp}}(\mathbb{S}^\iota_{M,i})} f_j^\xi$, for all $i\in\nn$, $((x_i)_{i=1}^\infty, (x_i^*)_{i=1}^\infty) $ is a $\varrho$-Schreier pair such that $(x_i)_{i=1}^\infty$ is a $1$-$\ell_1^{\varrho}$ spreading model and $(x^*_i)_{i=1}^\infty$ is a $1$-$c_0^{\varrho}$ spreading model.

\label{existence}
\end{theorem}

\begin{proof} We show first that there is $M\in[L]$ so that $(x_i)_{i=1}^\infty$ and $(x_i^*)_{i=1}^\infty$ satisfy the following.
\begin{enumerate}[$(i)$]

\item $(x_i)_{i=1}^\infty$ is a normalized block sequence in $X_\xi$ and $(x^*_i)_{i=1}^\infty$ is a normalized block sequence in $X_\xi^*$, 

\item $x^*_i(x_j)=0$ for all $i,j\in \nn$ with $i\neq j$, 

\item $x^*_i(x_i)=1$ for all $i\in\nn$, 

\item $\lim_i \|x_i\|_\beta=0$ for all $\beta<\iota$, 

\item $(x_i)_{i=1}^\infty$ is a $1$-$\ell_1^{\varrho}$ spreading model, 

\item $(x^*_i)_{i=1}^\infty$ is a $1$-$c_0^{\varrho}$ spreading model, 

\item $\text{\emph{supp}}(x_i), \text{\emph{supp}}(x^*_i)\in \mathcal{S}_\xi\cap \mathcal{S}_\iota$ for all $i\in\nn$. 

\end{enumerate}

Fix $N\in[L]$ such that $\mathcal{S}_\iota\cap[[L(1),\infty)]^{<\omega}\subset \mathcal{S}_\xi$.   By Proposition \ref{combinatorics}, there exists $K\in[\nn]$ such that $\cup_{i\in F}\text{supp}(\mathbb{S}^\iota_{N,K(i)})\in \mathcal{S}_\xi$ for any $F\in \mathcal{S}_{\varrho}$ and let $M=\cup_{i=1}^\infty \text{supp}(\mathbb{S}^\iota_{N,K(i)})\in [L]$.   Items $(i)$-$(iii)$ are clear. Item $(iv)$ follows from Corollary \ref{tail}.  Item $(v)$ follows from the fact that, with \[x_i=\sum_{j\in \text{supp}(\mathbb{S}^\iota_{M,i})} \mathbb{S}^\iota_{M,i}(j)e^\xi_j=\sum_{j\in \text{supp}(\mathbb{S}^\iota_{N,K(i)})}\mathbb{S}^\iota_{N,K(i)}(j)e^\xi_j,\] for any $F\in \mathcal{S}_{\varrho}$ and scalars $(a_i)_{i\in F}$  by Lemma \ref{repeatedaveragesfacts}, \[\Bigl\|\sum_{i\in F} a_i x_i\Bigr\|_{\ell_1}=\sum_{i\in F}|a_i|\] and \[\supp\Bigl(\sum_{i\in F}x_i\Bigr)\in \mathcal{S}_\xi,\] so  \[\Bigl\|\sum_{i\in F}a_ix_i\Bigr\|_\xi=\sum_{i\in F}|a_i|.\] We therefore deduce that $(x_i)_{i=1}^\infty$ is a $1$-$\ell_1^{\varrho}$ spreading model.

Item $(vi)$ follows from the fact that, with \[x_i^*=\sum_{j\in \text{supp}(\mathbb{S}^\iota_{M,i})}f^\xi_j=\sum_{j\in \text{supp}(\mathbb{S}^\iota_{N,K(i)})}f^\xi_j,\] for any $F\in \mathcal{S}_{\varrho}$, \[\Bigl\|\sum_{i\in F} x^*_i\Bigr\|_{c_0}=1\] and \[\supp\Bigl(\sum_{i\in F}x^*_i\Bigr)\in \mathcal{S}_\xi,\] \[\Bigl\|\sum_{i\in F}x^*_i\Bigr\|_\xi= 1.\]    By $1$-unconditionality, $(x^*_i)_{i=1}^\infty$ is a $1$-$c_0^{\varrho}$ spreading model. 

Item $(vii)$ follows from the fact that $\supp(\mathbb{S}^\iota_{R,j})\in \mathcal{S}_\iota$ for all $R\in[\nn]$ and $j\in\nn$, together with the fact that $\mathcal{S}_\iota\cap [[L(1), \infty)]^{<\omega}\subset \mathcal{S}_\xi$ and $M\in [L]$.    

Now we apply Lemma \ref{lastlegs1} to finish the proof.

\end{proof}

\begin{corollary}
Fix $\xi<\omega_1$ and  $\varrho\in R(\xi)$.  Suppose that  $(x_i)_{i=1}^\infty\subset X_\xi$ is a normalized block sequence which is equivalent to a subsequence of $(e^\varrho_i)_{i=1}^\infty$. Then for any $\delta\in R(\xi)\cap \varrho+$, there exists a convex block sequence $(y_i)_{i=1}^\infty$ of $(x_i)_{i=1}^\infty$ such that $(y_i)_{i=1}^\infty$ is equivalent to a subsequence of $(e^\delta_i)_{i=1}^\infty$. 

\end{corollary}

\begin{proof}
Since $(x_i)_{i=1}^\infty$ is equivalent to a subsequence of $(e^\varrho_i)_{i=1}^\infty$, it is sufficient to show that any $\delta\in R(\xi)\cap \varrho+$,  any subsequence of $(e^\varrho_i)_{i=1}^\infty$ has a convex block sequence which is equivalent to a subsequence of $(e^\delta_i)_{i=1}^\infty$. Since $R(\xi)\cap \varrho+=R(\varrho)$, this follows from Theorem \ref{existence}. 

\end{proof}

In the next theorem, we say a sequence $(y^*_i)_{i=1}^\infty$ is an $\ell_\infty^+$ \emph{block sequence} of $(x^*_i)_{i=1}^\infty$ if there exist sets $E_1<E_2<\ldots$ such that $y^*_i=\sum_{j\in E_i}x^*_j$.

\begin{corollary} Fix $\xi<\omega_1$ and $\varrho\in R(\xi)$. Suppose that $(x_i^*)_{i=1}^\infty\subset X_\xi^*$ is a normalized block sequence which is equivalent to a subsequence of $(f^\varrho_i)_{i=1}^\infty$. Then for any $\delta\in R(\xi)\cap \varrho+=R(\varrho)$, there exists an $\ell_\infty^+$ block sequence $(y_i^*)_{i=1}^\infty$ of $(x_i^*)_{i=1}^\infty$ such that $(y_i^*)_{i=1}^\infty$ is equivalent to a subsequence of $(f^\delta_i)_{i=1}^\infty$.  

\end{corollary}

\begin{proof} By Corollary \ref{moveup2}, by passing to a subsequence and relabeling, we can assume that there exists a normalized block sequence  $(x_i)_{i=1}^\infty\subset X_\xi$  such that $((x_i)_{i=1}^\infty, (x_i^*)_{i=1}^\infty)$ is a $\varrho$-Schreier pair. 

Since $\delta=R(\xi)\cap \varrho+=R(\varrho)$, by Theorem \ref{existence}, there exist $I\in [\nn]$,  $E_1<E_2<\ldots$, $E_i\in [J]^{<\omega}$, and $(z_i)_{i=1}^\infty\in \prod_{i=1}^\infty \text{co}\{e^\varrho_j:j\in E_i\}$ such that, with $z^*_i=\sum_{j\in E_i}f^\varrho_j$, 

\begin{enumerate}[$(i)$]\item $(z_i)_{i=1}^\infty\approx (e^\delta_i)_{i\in I}$, 

\item $(z^*_i)_{i=1}^\infty\approx (e^\delta_i)_{i\in I}$.  

\end{enumerate}

Then if $T=\sum_{i=1}^\infty x_i\otimes f^\varrho_{J(i)}$ and $S=\sum_{i=1}^\infty x^*_i\otimes e^\varrho_{J(i)}$, $(y_i)_{i=1}^\infty:=(Tz_i)_{i=1}^\infty\approx (e^\delta_i)_{i\in I}$ is a convex block sequence of $(x_i)_{i=1}^\infty$ and $(y^*_i)_{i=1}^\infty:=(Sz^*_i)_{i=1}^\infty\approx (f^\delta_i)_{i\in I}$ is an $\ell_\infty^+$ block sequence of $(x^*_i)_{i=1}^\infty$.  

\end{proof}

\section{Strictly singular operators on Schreier spaces}

In this section we apply the results classifying block sequences in Schreier spaces, mainly Theorem \ref{final1} and Corollary \ref{existence}, in the study the classes of strict singular operators on Schreier spaces. 

We recall first different classes of strictly singular operators. Let $X,Y$ be Banach spaces. An operator $T\in\mathcal{L}(X, Y)$ is strictly singular if its restriction to any infinite dimensional subspace of $X$ is not an isomorphism onto its image, in other words any infinite dimensional subspace $Z\subset X$ contains for any $\varepsilon>0$ a vector $x$ with $\|x\|=1$ and $\|Tx\|<\varepsilon$. One can quantify this fact using countable ordinals; following \cite{ADST} we say that an operator $T\in\mathcal{L}(X,Y)$ is  $\mathcal{S}_\xi$-strictly singular, $\xi<\omega_1$, if for any $\varepsilon>0$ and any basic sequence $(x_i)_{i=1}^\infty\subset X$ there is $F\in\mathcal{S}_\xi$ and $x\in \mathrm{span}\{x_i: i\in F\}$ with $\|x\|=1$ and $\|Tx\|<\varepsilon$. Finally, the most restrictive form of strict singularity appearing in the literature "localizes" the choice of vectors testing strict singularity; we say that an operator $T: X\to Y$ is finitely strictly singular (or super strictly singular, cf. \cite{FHR,P}), if for any $\varepsilon>0$ there is $n\in\nn$ such that for any $n$-dimensional subspace $E$ of $X$ there is $x\in E$ with $\|x\|=1$ and $\|Tx\|<\varepsilon$. 

With the above notation $\mathcal{S}_0$-strictly singular operators are precisely the compact ones. The class of strictly singular operators in $\mathcal{L}(X,Y)$ is denoted by $\mathcal{SS}(X,Y)$, the class of $\mathcal{S}_\xi$-strictly singular ones, $\xi<\omega_1$ - by $\mathcal{SS}_\xi(X,Y)$, the class of finitely strictly singular operators - by $\mathcal{FSS}(X,Y)$ and finally, the class of compact operators - by $\mathcal{K}(X_\xi)$. For $Y=X$ we skip $Y$ in the notation and write $\mathcal{SS}(X)$ etc. Families of strictly singular, finitely strictly singular and compact operators form closed operator ideals (i.e.  closed in the norm topology, two-sided ideals which are also vector spaces), contrary to the case of $\xi$-strictly singular operators, see \cite{ADST}. We shall use the following useful properties.

\begin{proposition}\cite{ADST} \label{strictlysingularfacts}    
Let $X,Y$ be Banach spaces. Then for any $0<\zeta\leqslant\xi<\omega_1$, \[\mathcal{K}(X,Y)\subset\mathcal{FSS}(X,Y)\subset\mathcal{SS}_\zeta(X,Y)\subset\mathcal{SS}_\xi(X,Y)\subset\mathcal{SS}(X,Y).\] 
\end{proposition}

\begin{lemma}\label{ss reduced to block sequences} Let $\xi<\omega_1$ and  $T\in\mathcal{L}(X,Y)$, where $X$ is a Banach space with a basis, not containing an isomorphic copy of $\ell_1$. Then $T\in\mathcal{SS}_\xi(X,Y)$ iff for any $\varepsilon>0$ and any block sequence $(x_i)_{i=1}^\infty\subset X$ there is $F\in\mathcal{S}_\xi$ and $x\in \mathrm{span}\{x_i: i\in F\}$ with $\|x\|=1$ and $\|Tx\|<\varepsilon$.    
\end{lemma}
\begin{proof}
For $\xi=0$ (i.e. for compact operators) use Odell's theorem \cite[p. 377]{R} stating that an operator on a Banach space not containing $\ell_1$ is compact iff is completely continuous, i.e. carries weakly null sequences to norm null sequences. 

Fix $0<\xi<\omega_1$ and pick any normalized basic sequence $(x_i)_{i=1}^\infty\subset X$. By Rosenthal $\ell_1$ theorem, passing to a subsequence and allowing some perturbation we can assume that $(x_{2i}-x_{2i+1})_{i=1}^\infty$ is a block sequence. By the assumption on $T$ pick $F\in\mathcal{S}_\xi$ and $x\in\mathrm{span}\{x_{2i}-x_{2i+1}: i\in F\}$ with $\|x\|=1$ and $\|Tx\|<\varepsilon$. To finish the proof we use the following observation: for any $F\in\mathcal{S}_\xi$, $\{2i, 2i+1: i\in F\}\in \mathcal{S}_\xi$, which can be proved easily by induction on $0<\xi<\omega_1$. 
\end{proof}

We start with the following result on strictly singular operators between Schreier spaces and arbitrary Banach spaces and later proceed to the case of operators between Schreier spaces. 
\begin{proposition} Let $Y$ be a Banach space.  Then 
\begin{enumerate}
    \item 
$\mathcal{SS}(X_1,Y)=\mathcal{FSS}(X_1,Y)$,
    \item $\mathcal{SS}(X_\xi,Y)=\mathcal{SS}_\xi(X_\xi,Y)$, for any $\xi<\omega_1$.
\end{enumerate}
\end{proposition}
\begin{proof} 
(1) Fix $T\in\mathcal{SS}(X_1,Y)$. Take any finite dimensional subspaces $(E_i)_{i=1}^\infty$ of $X_1$ with $\mathrm{dim}(E_i)\to\infty$.
By \cite{M}, see also \cite[Lemma 13]{CP} in any $k$-dimensional subspace of $X_1$ there is a vector $x$ which is $k$-flat, i.e. such that $\|x\|_0=1$ and at least $k$ many coordinates equal one. Therefore we can choose $(k_i)_{i=1}^\infty\subset\nn$ and a block sequence $(x_i)_{i=1}^\infty\subset X_1$ such that $x_i\in E_{k_i}$ and $x_i$ is $i$-flat for each $i\in\nn$. Assume that $(Tx_i)_{i=1}^\infty$ is seminormalized, otherwise the proof is finished. As $(x_i)_{i=1}^\infty$ is weakly null, again  passing to a subsequence we can assume that $(Tx_i)_{i=1}^\infty$ is a basic sequence. Let $y_i=\frac{x_i}{\|x_i\|_1}\in E_{k_i}$ for any $i\in\nn$. As  $\|x_i\|_1\geq i$, $\|y_i\|_0\to 0$, thus $(y_i)_{i=1}^\infty$ contains a subsequence equivalent to the unit vector basis of $c_0$ (\cite[Prop. 4.3]{GL}, see also Theorem \ref{final1}). Thus $\liminf_{i\to\infty}\|Ty_i\|=0$ by strict singularity  of $T$. As subspaces $(E_i)_{i=1}^\infty$ were  chosen arbitrary, we proved the finite strict singularity of $T$.

(2) Fix $T\in\mathcal{SS}(X_\xi,Y)$. Take any normalized block sequence $(x_i)_{i=1}^\infty\subset X_\xi$. 
Assume that $(Tx_i)_{i=1}^\infty$ is seminormalized, otherwise the proof is finished. As $(x_i)_{i=1}^\infty$ is weakly null, again  passing to a subsequence we can assume that $(Tx_i)_{i=1}^\infty$ is a basic sequence. By Theorem \ref{final1},  passing to a further subsequence we can assume also that $(x_i)_{i=1}^\infty$ is equivalent to a subsequence $(e_i^\varrho)_{i\in M}$ of the unit vector basis of $X_\varrho$ for some $\varrho\in R(\xi)$. We can also assume that $\mathcal{S}_\varrho\cap [M]^{<\infty}\subset \mathcal{S}_\xi$ by Proposition \ref{schreierfamilyfacts}$(iv)$. Applying  Corollary \ref{isometricc0}, passing to a further subsequence of $M$ if necessary, we obtain a sequence $(y_n)_{n=1}^\infty=(\sum_{i=1}^\infty\mathbb{S}^\varrho_{M,n}(i)x_i)_{n=1}^\infty$ that is equivalent to the unit vector basis of $c_0$. By strict singularity of $T$, $\|Ty_n\|\to 0$. As $\mathrm{supp} (\mathbb{S}_{M,n}^\varrho)\in\mathcal{S}_\varrho$, $n\in\nn$ (Lemma \ref{repeatedaveragesfacts}$(ii)$), by the choice on $M$,  $\mathrm{supp} (\mathbb{S}_{M,n}^\varrho)\in\mathcal{S}_\xi$, $n\in\nn$, which ends the proof.  
\end{proof}

For characterization of the order of strict singular operators on Schreier spaces we shall use the following lemma.

\begin{lemma}
Let $\xi,\zeta<\omega_1$ and $T\in\mathcal{L}(X_\xi,X_\zeta)$. 
\begin{enumerate}
\item[(i)] For any normalized block sequence $(x_i)_{i=1}^\infty\subset X_\xi$ with $(Tx_i)_{i=1}^\infty$ seminormalized there are $\mu\in R(\xi)$, $\varrho\in R(\zeta)$, and $M\in[\nn]$ such that $(x_i)_{i\in M}$ is equivalent to $(e_i^\mu)_{i\in M}$ and $(Tx_i)_{i\in M}$ is equivalent to $(e_i^\varrho)_{i\in M}$.
\item[(ii)] $T\in\mathcal{SS}(X_\xi,X_\zeta)$ iff for any block sequence $(x_i)_{i=1}^\infty\subset X_\xi$ equivalent to a subsequence of $(e_i^\mu)_{i=1}^\infty$ with $(Tx_i)_{i=1}^\infty$  equivalent to a subsequence of $(e_i^\varrho)_{i=1}^\infty$ we have $\varrho<\mu$.
\item[(iii)] Assume $T\in\mathcal{SS}(X_\xi,X_\zeta)\setminus\mathcal{K}(X_\xi,X_\zeta)$ and let $\eta\in R(\zeta)$ be maximal with the following property:  for some normalized block sequence $(x_i)_{i=1}^\infty\subset X_\xi$, $(Tx_i)_{i=1}^\infty$ is equivalent to a subsequence of $(e_i^{\eta})_{i=1}^\infty$. Then    
$T\in\mathcal{SS}_{\eta+1}(X_\xi, X_\zeta)\setminus\mathcal{SS}_{\eta}(X_\xi, X_\zeta)$.
\end{enumerate}

\label{techlemmaforoperators}
\end{lemma}
\begin{proof} To prove $(i)$, use Theorem \ref{final1} to pick $\mu\in R(\xi)$, $N=(n_i)_{i=1}^\infty,L=(l_i)_{i=1}^\infty\in[\nn]$  with $(x_m)_{m\in M}$ equivalent to $(e_l^\mu)_{l\in L}$ and $\varrho\in R(\zeta)$, $M=(m_i)_{i=1}^\infty\in [M]$, $J=(j_i)_{i=1}^\infty\in[\nn]$ with $(Tx_{m})_{m\in M}$ equivalent to $(e_{j}^\varrho)_{j\in J}$. Passing to a further subsequence we can assume that $\max\{m_i,l_i,j_i \}<\min\{m_{i+1}, l_{i+1}, j_{i+1}\}$ for all $i\in\nn$. By Proposition \ref{schreierspacefacts}$(iii)$ $M$ satisfies the assertion of $(i)$.

We proceed to the proof of $(ii)$. Assume $T$ is strictly singular and pick $(x_i)_{i=1}^\infty$, $\mu\in R(\xi)$ and $\varrho\in R(\zeta)$ as in $(ii)$. Apply first $(i)$ obtaining $M\in[\nn]$ with $(x_i)_{i\in M}$ equivalent to $(e_i^\mu)_{i\in M}$ and $(Tx_i)_{i\in M}$ equivalent to $(e_i^\varrho)_{i\in M}$. If $\varrho>\mu$, then Prop. \ref{smallbetanorm} would yield contradiction with boundedness of $T$, whereas strict singularity of $T$ ensures $\varrho\neq\mu$. 

Assume $\varrho>\mu$ for any $\varrho$ and $\mu$ as in $(ii)$. Pick any normalized block sequence $(x_i)_{i=1}^\infty\subset X_\xi$. If $\liminf_i\|Tx_i\|=0$ the proof is finished, otherwise apply $(i)$ obtaining $\varrho\in R(\zeta)$, $\mu\in R(\xi)$ and $M\in[\nn]$ with $(x_i)_{i=1}^\infty$ $C$-equivalent to $(e^\mu_i)_{i\in M}$ and $(Tx_i)_{i\in M}$ $C$-equivalent to $(e^\varrho_i)_{i\in M}$ for some $C>0$. By the assumption $\mu<\varrho$, thus for any $\varepsilon>0$ by Proposition \ref{smallbetanorm} there is $u=(a_i)_{i=1}^\infty\in c_{00}$ with $\supp(u)\in \mathcal{S}_\mu\cap [M]^{<\infty}$, $\sum_{i=1}^\infty|a_i|=1$ and $\|u\|_\varrho<\varepsilon$. Then 
\[\Big\|\sum_{i\in M} a_ix_i\Big\|_\xi\geqslant C^{-1}\Big\|\sum_{i\in M}a_ie_i^\mu\Big\|_\mu=C^{-1}\]
whereas
\[\Big\|\sum_{i\in M}a_iTx_i\Big\|_\zeta\leqslant C\Big\|\sum_{i\in M}a_ie^\varrho_i\Big\|_\varrho<C\varepsilon\]
As $\varepsilon>0$ is arbitrary, by Lemma \ref{ss reduced to block sequences} we finish the proof of strict singularity of $T$. 

To show $(iii)$ pick $\eta\in R(\xi)$ as in $(iii)$ and note $\eta<\xi$ by $(ii)$. First pick a normalized block sequence $(x_i^0)_{i=1}^\infty\subset X_\xi$ with $(Tx_i^0)_{i=1}^\infty$ $C$-equivalent to a subsequence of $(e^\eta_i)_{i=1}^\infty$ for some $C>0$. By continuity of $T$, up to some perturbation, we can assume that $(Tx_i^0)_{i=1}^\infty$ is also a block sequence. Note that for any $E\in\mathcal{S}_\eta$ and scalars $(a_i)_{i\in E}$ we have
\[\Big\|T\sum_{i\in E}a_ix_i^0\Big\|_\zeta\geqslant C^{-1}\Big\|\sum_{i\in E}a_ie_i^\eta\Big\|_\eta= C^{-1}\sum_{i\in E}|a_i|\geqslant C^{-1}\Big\|\sum_{i\in E}a_ix_i^0\Big\|_\xi\]
thus $T\not\in\mathcal{SS}_\eta(X_\xi,X_\zeta)$. 

We prove now that $T\in\mathcal{SS}_{\eta+1}(X_\xi,X_\zeta)$. Take any normalized block sequence $(x_i)_{i=1}^\infty\subset X_\xi$. We repeat slightly changed proof of $(ii)$. If $\liminf_i\|Tx_i\|=0$ the proof is finished, otherwise apply $(i)$ obtaining $\varrho\in R(\zeta)$, $\mu\in R(\xi)$ and $M\in[\nn]$ with $(x_i)_{i=1}^\infty$ $C$-equivalent to $(e^\mu_i)_{i\in M}$ and $(Tx_i)_{i\in M}$ $C$-equivalent to $(e^\varrho_i)_{i\in M}$ for some $C>0$. By the assumption $T$ is strictly singular, thus by $(ii)$ $\varrho<\mu$. Therefore we can assume by Proposition \ref{schreierfamilyfacts}$(iv)$ and the definition of $\eta$ that $\mathcal{S}_{\varrho+1}\cap [M]^{<\infty}\subset \mathcal{S}_\mu\cap \mathcal{S}_{\eta+1}$. Now for any $\varepsilon>0$ by Proposition \ref{smallbetanorm} pick $u=(a_i)_{i=1}^\infty\in c_{00}$ with $\supp(u)\in \mathcal{S}_{\varrho+1}\cap [M]^{<\infty}$, $\sum_{i=1}^\infty|a_i|=1$ and $\|u\|_\varrho<\varepsilon$. Then 
\[\Big\|\sum_{i\in M} a_ix_i\Big\|_\xi\geqslant C^{-1}\Big\|\sum_{i\in M}a_ie_i^\mu\Big\|_\mu=C^{-1}\]
whereas
\[\Big\|\sum_{i\in M}a_iTx_i\Big\|_\zeta\leqslant C\Big\|\sum_{i\in M}a_ie^\varrho_i\Big\|_\varrho<C\varepsilon\]
As $\varepsilon>0$ is arbitrary and $\supp(u)\in\mathcal{S}_{\eta+1}$ we finish the proof of $(\eta+1)$-strict singularity of $T$. 
\end{proof}

The following observation provides, for any $\zeta<\xi<\omega_1$, operators that are  $\mathcal{S}_\xi$-strictly singular, but not  $\mathcal{S}_\zeta$-strictly singular, completing the range of examples presented in \cite[Example 2.7]{ADST}.

\begin{proposition}    
Fix $\zeta\leqslant\xi<\omega_1$. Let $id_{\xi,\zeta}$ be the formal identity mapping $id_{\xi,\zeta}: X_\xi\to X_\zeta$. Then the following are equivalent
\begin{enumerate}
\item $\zeta<\xi$,
\item $id_{\xi,\zeta}\in\mathcal{SS}_{\zeta+1}(X_\xi,X_\zeta)\setminus\mathcal{SS}_\zeta(X_\xi,X_\zeta)\subset\mathcal{SS}_\xi(X_\xi,X_\zeta)\setminus\mathcal{SS}_\zeta(X_\xi,X_\zeta)$. 
\end{enumerate}
\label{identityschreierss}
\end{proposition}

\begin{proof} Note that (2) implies strict singularity of $id_{\xi,\zeta}$, thus also trivially (1). 

Assume now (1). In order to show (2) it is enough to prove, by Lemma \ref{techlemmaforoperators}$(iii)$ and Proposition \ref{strictlysingularfacts}, that $id_{\xi,\zeta}$ is strictly singular. We shall prove it applying Lemma \ref{techlemmaforoperators}$(ii)$. Pick any normalized block sequence $(x_i)_{i=1}^\infty\subset X_\xi$ equivalent to a subsequence of $(e_i^\mu)_{i=1}^\infty$ with $(x_i)_{i=1}^\infty\subset X_\zeta$ equivalent to a subsequence of $(e^\varrho_i)_{i=1}^\infty$ for some $\mu\in R(\xi)$ and $\varrho\in R(\zeta)$ (in particular $(x_i)_{i=1}^\infty$ is seminormalized in $X_\zeta$). Then by Lemma \ref{lastlegs1} $\mu=\xi-\iota$ and $\varrho=\zeta-\iota$, where $\iota=\min \{\beta\leqslant \zeta: \liminf_i\|x_i\|_\beta>0\}$. As $\zeta<\xi$ by assumption, also $\varrho<\mu$, which by Lemma \ref{ss reduced to block sequences} ends the proof of strict singularity of $id_{\xi,\zeta}$ and thus of (2). 
\end{proof}

The following main result of this section can be regarded as a non-hereditary version of results of \cite{ABM} for families of strictly singular operators of higher order.

\begin{theorem}
Fix $1\leqslant\xi<\omega_1$ and let $R(\xi)=\{\varrho_0,\dots,\varrho_k\}$ with $0=\varrho_0<\dots<\varrho_k=\xi$, Then 
\begin{enumerate}
\item \[\mathcal{SS}(X_\xi)\setminus\mathcal{K}(X_\xi)=\bigcup_{s=0}^{k-1}\Big(\mathcal{SS}_{\varrho_s+1}(X_\xi)\setminus\mathcal{SS}_{\varrho_s}(X_\xi)\Big)\]
with $\mathcal{SS}_{\varrho_s+1}(X_\xi)\setminus\mathcal{SS}_{\varrho_s}(X_\xi)\neq\emptyset$ for each $0\leqslant s<k$.
\item 
The family $\mathcal{SS}_\varrho(X_\xi)$ is a closed operator ideal for any $\varrho<\omega_1$. 
\item The product of any $(k+1)$ strictly singular operators in $\mathcal{L}(X_\xi)$ is compact. 
\item There are strictly singular operators $T_1,\dots,T_k\in\mathcal{L}(X_\xi)$ with a non-compact product $T_k\cdots T_1$. 
\end{enumerate}
\label{structure of ss operators}
\end{theorem}
\begin{proof} 
(1) The equality follows immediately from Lemma \ref{techlemmaforoperators}$(iii)$. Fix $0\leqslant s<k$. By Theorem \ref{existence} there is a normalized sequence $(x_i)_{i=1}^\infty\subset X_\xi$ equivalent to some subsequence of $(e_i^{\varrho_s})_{i\in\nn}$. By Prop. \ref{identityschreierss} the operator $T\in\mathcal{L}(X_\xi)$ carrying each $e_i^\xi$ to $x_i$, $i\in\nn$, is in the class $\mathcal{SS}_{\varrho_s+1}(X_\xi)\setminus \mathcal{SS}_{\varrho_s}(X_\xi)$.

(2) By (1) for any $\varrho<\omega_1$, $\mathcal{SS}_\varrho(X_\xi)=\mathcal{SS}_{\varrho_s+1}(X_\xi)$, where $\varrho_s$ is maximal in $R(\xi)\setminus \{\xi\}$ with $\varrho_s<\varrho$. 

In order to prove that each family  $\mathcal{SS}_{\varrho_s+1}(X_\xi)$, $0\leqslant s<k$, is a closed operator ideal it is enough to show that is is closed under addition, as all the other properties were shown in \cite[Prop. 2.4]{ADST}. Pick $T,S\in\mathcal{SS}_{\varrho_s+1}(X_\xi)$. 
Let $(x_i)_{i=1}^\infty\subset X_\xi$ be a normalized block sequence with $((T+S)(x_i))_{i=1}^\infty$ equivalent to some subsequence of $(e_i^\varrho)_{i=1}^\infty$ with $\varrho\in R(\xi)$. We shall prove that $\varrho\leqslant\varrho_s$, which by Lemma \ref{techlemmaforoperators}$(iii)$ yields $T+S\in\mathcal{SS}_{\varrho_s+1}(X_\xi)$. If either of $(Tx_i)_{i=1}^\infty$ or  $(Sx_i)_{i=1}^\infty$ are not seminormalized, say $\liminf_i\|Tx_i\|_\xi=0$,  passing to a subsequence we can assume that $((T+S)x_i)_{i=1}^\infty$ is equivalent to $(Sx_i)_{i=1}^\infty$. By Theorem \ref{final1} some subsequence of $(Sx_i)_{i=1}^\infty$ is equivalent to a subsequence of $(e^\mu_i)_{i=1}^\infty$ for some $\mu\in R(\xi)$. By Lemma \ref{techlemmaforoperators}$(iii)$, $\mu\leqslant\varrho_s$, whereas $\varrho=\mu$ by Proposition \ref{smallbetanorm}. Assume that both $(Tx_i)_{i=1}^\infty$ and $(Sx_i)_{i=1}^\infty$ are seminormalized. By Theorem \ref{final1}, passing to a further subsequence we can assume that $(Tx_i)_{i=1}^\infty$ is equivalent to a subsequence of $(e_i^\eta)_{i=1}^\infty$  and $(Sx_i)_{i=1}^\infty$ is equivalent to a subsequence of $(e_i^\mu)_{i=1}^\infty$ for some $\eta,\mu\in R(\xi)$. By Lemma \ref{techlemmaforoperators}$(iii)$, $\eta,\mu\leqslant\varrho_s$. As the bases of Schreier spaces are 1-right dominant (Proposition \ref{schreierspacefacts}$(i)$) some subsequence of $(e_i^{\varrho_s})_{i=1}^\infty$ dominates both $(Tx_i)_{i=1}^\infty$ and $(Sx_i)_{i=1}^\infty$, thus also $((T+s)x_i)_{i=1}^\infty$ and hence a subsequence of $(e_i^\varrho)_{i=1}^\infty$. Again Proposition \ref{smallbetanorm} yields $\varrho\leqslant\varrho_s$.

(3) We follow the approach of \cite[Thm 4.1]{ADST}. Pick $T_1,\dots,T_{k+1}\in\mathcal{SS}(X_\xi)$. By Lemma \ref{ss reduced to block sequences} it is enough to show that for any normalized block sequence $(x_i)_{i=1}^\infty\subset X_\xi$, $\liminf_i\|T_{k+1}\cdots T_1x_i\|=0$. Assume the opposite. Passing to a subsequence, allowing perturbation and applying  Theorem \ref{final1}, we can assume that  each sequence  $(T_s\cdots T_1x_i)_{i=1}^\infty$, $s=1,\dots,k+1$, is a block sequence equivalent to a subsequence of $(e_i^{\mu_s})_{i=1}^\infty$, for some $\mu_s\in R(\xi)$. By Lemma \ref{techlemmaforoperators}$(ii)$,  $\mu_1<\dots<\mu_{k+1}<\xi$ which yields a contradiction.  

(4) For any $0\leqslant s\leqslant k$ by Theorem \ref{existence} pick a normalized block sequence $(x_i^s)_{i=1}^\infty\subset X_\xi$ equivalent to a subsequence  $(e_j^{\varrho_s})_{j\in J_s}$ of the unit vector basis of $X_{\varrho_s}$ and spanning a complemented subspace of $X_\xi$ (by the definition of a $\varrho$-Schreier pair). Let $J_s=(j_i^s)_{i=1}^\infty$, $0\leqslant s\leqslant k$. Passing to subsequences we can assume that $\max\{j_i^0,\dots,j_i^{k-1}\}<\min \{j_{i+1}^0,\dots,j_{i+1}^{k-1}\}$ for any  $i\in\nn$, thus by Proposition \ref{schreierfamilyfacts}$(iii)$ each sequence $(x_i^s)_{i=1}^\infty$ is equivalent to $(e_j^{\varrho_s})_{j\in J_0}$, $0\leqslant s\leqslant k$. By Propositions \ref{schreierfamilyfacts}$(iv)$ and \ref{identityschreierss} for any $1\leqslant s\leqslant k$ the mapping carrying $x_i^s$ to $x_i^{s-1}$, $i\in\nn$, extends to a bounded strictly singular operator $\overline{\mathrm{span}}\{x_i^s: i\in\nn\}\to \overline{\mathrm{span}}\{x_i^{s-1}: i\in\nn\}$, call it $R_s$. Finally let $T_s:=R_sP_s\in\mathcal{SS}(X_\xi)$, where $P_s\in\mathcal{L}(X_\xi)$ is the projection onto $\overline{\mathrm{span}}\{x_i^s: i\in\nn\}$, for any $1\leqslant s\leqslant k$. Note that $T_k\cdots T_1$ is non-compact as $T_k\cdots T_1x_i^k=x_i^0$ for any $i\in\nn$, which ends the proof. 
\end{proof}

\section{Operator ideals on Schreier spaces, their duals and biduals - examples}

Recently the maximal cardinality, i.e. $2^\mathfrak{c}$, of the lattice of closed operator ideals was shown for a list of classical separable Banach spaces. Among the spaces with this property are in particular $\mathcal{L}(L_p)$, $1<p<\infty$, $p\neq 2$ \cite{JS}, $\ell_p\oplus\ell_q$, $1\leqslant p,q\leqslant\infty$, apart from $(p,q)=(1,\infty)$ \cite{FSZ}, Schlumprecht space and Schreier spaces of finite order \cite{MP,Schl}, 
Schreier spaces of arbitrary order \cite{Schl}. In this section we extend the method of \cite{MP} obtaining $2^\mathfrak{c}$ many small (i.e. consisting of strictly singular operators) closed operator ideals on Schreier spaces of arbitrary order, their dual and bidual spaces. We follow the reasoning presented in \cite{MP} for Schreier spaces of finite order, substituting the crucial step based on \cite[Thm 1.1]{GL} by an analogous combinatorial reasoning using the notion of $\xi$-injective map. This combinatorial part,  relying heavily on Lemmas \ref{mod1} and \ref{weaksumming}, follows part of the scheme of \cite{GL}, does not provide generalization of \cite[Thm 1.1]{GL} to the case of Schreier spaces of arbitrary order, but is sufficient for our purposes. 


We start with the generalization of some results of \cite{GL}. We fix the following notation. 
\begin{definition}
For any $\xi<\omega_1$ and finite $A\subset\mathbb{N}$ let
\[\tau_\xi(A)=\min\{k\in\nn: \text{ there are successive } F_1,\dots,F_k\in\mathcal{S}_\xi \text{ with } A=F_1\cup\dots\cup F_k\}.\]
\end{definition}
As the above is a simple reformulation of the parameter $\tau_\xi$ introduced in \cite{GL} we keep the same notation. The next fact, based on Lemma \ref{mod1}, generalizes \cite[Lemma 3.2 (5)]{GL}, shown for  $\xi<\omega$ with a constant depending on $\xi$.
\begin{lemma}
Let $A=\bigcup_{i=1}^nA_i\in[\nn]^{<\infty}$. Then for any $\xi<\omega_1$ we have
\[\tau_\xi(A)\leqslant \sum_{i=1}^n\tau_\xi(A_i)\]  
\label{subadditivity}
\end{lemma}
\begin{proof} We can assume that $(A_i)_{i=1}^n$ are pairwise disjoint sets. 
For any $i=1,\dots,n$ let $k_i=\tau_\xi(A_i)$ and pick $(F_{i,k})_{k=1}^{k_i}\subset\mathcal{S}_\xi$ with $A_i=\bigcup_{i=1}^{k_i}F_{i,k}$. By  the property  $P_\xi$ in the proof of Lemma \ref{mod1}, applied to the family $\{F_{i,k}: k=1,\dots,k_n, i=1,\dots,n\}$, there is a sequence of successive $\mathcal{S}_\xi$-sets $(G_k)_{k=1}^{k_1+\dots+k_n}$ with $G_1\cup\dots\cup G_{k_1+\dots+k_n}=A_1\cup\dots\cup A_n=A$, which ends the proof. 
\end{proof}

The following notion of $\xi$-injective maps is crucial in our proof. Maps of this type appear in \cite[Prop. 3.11]{GL}, and roughly speaking, in the case of $\xi<\omega$, can be reduced to canonical increasing bijections between infinite sets of integers (cf. \cite[Prop. 3.12]{GL}). \begin{definition} Fix $L,M\in[\nn]$ and $\xi<\omega_1$.
We say that a map $\psi: L\to M$ is $\xi$-injective if 
\[
\sup_{F\in\mathcal{S}_\xi}\tau_\xi (\psi^{-1}(F))<\infty
\]
\end{definition}

\begin{lemma} Fix $\xi<\omega_1$,  $L,M\in[\nn]$ and a map  $\psi: L\to M$. Then the following are equivalent.
\begin{enumerate}
    \item $\psi$ is $\xi$-injective,
    \item there is a constant $C>0$ such that for any $F\in[M]^{<\infty}$ we have
$\tau_\xi(\psi^{-1}(F))\leqslant C\tau_\xi(F)$.
\end{enumerate} 
\label{lemmaxiinjective}
\end{lemma}
\begin{proof}
The implication $(2)\Rightarrow (1)$ is obvious. For $(1)\Rightarrow (2)$  let $C=\sup_{F\in\mathcal{S}_\xi}\tau_\xi(\psi^{-1}(F))$. Given  $F\in[M]^{<\infty}$, let $k=\tau_\xi(F)$ and pick successive $\mathcal{S}_\xi$ sets $F_1<\dots<F_k$ with $F=\bigcup_{i=1}^kF_i$. Then by Lemma \ref{subadditivity}
\[
\tau_\xi(\psi^{-1}(F)) \leqslant \sum_{i=1}^k\tau_\xi(\psi^{-1}(F_i))\leqslant kC=C\tau_\xi(F).    
\]   
\end{proof}

The next result, based on Lemma \ref{weaksumming}, partly generalizes  \cite[Prop. 3.11]{GL}; we present it below for any $\xi<\omega_1$, but with $\xi$-injectivity of $\psi$ only instead of $\zeta$-injectivity of $\psi$ for all $\zeta\leqslant\xi$.
\begin{proposition}
Fix $L,M\in[\nn]$ and $\xi<\omega_1$. Assume there is a map $\psi: L\to M$ and a bounded operator $R: X_\xi(L)\to X_\xi(M)$ with $Re_l=Re_{\psi(l)}$ for any $l\in L$. Then $\psi$  is $\xi$-injective.      
\label{reductionoptomap}
 \end{proposition}
\begin{proof} We present a version of reasoning in \cite{GL}. Assume that $\psi$ is not $\xi$-injective, thus there is a sequence of successive $\mathcal{S}_\xi$   sets $(F_n)_{n=1}^\infty$ such that $(\psi^{-1}(F_n))_{n=1}^\infty$ is a sequence of successive sets and  $\tau_\xi(\psi^{-1}(F_n))>n$ for each $n\in\nn$. For any $n\in\nn$ denote by $(E_i)_{i\in I_n}$, $|I_n|=n$, the sequence of successive maximal $\mathcal{S}_\xi$ sets with $E_i\subset F_n$ for all $i\in I_n$. We assure also that $\bigcup_{n=1}^\infty I_n=\nn$. 

Let $M=\bigcup_{n=1}^\infty\bigcup_{i\in I_n}E_i$, then $F_i^\xi(M)=E_i$ for all $i\in\nn$. Take the sequence of repeated averages $(\mathbb{S}_{M,i}^\xi)_{i=1}^\infty$ and note that $\supp (\mathbb{S}_{M,i}^\xi)=E_i$ for all $i\in\nn$. Let $x_i=\sum_{j\in E_i}\mathbb{S}_{M,i}^\xi(j)e_j^\xi$ for all $i\in\nn$ and $y_n=\frac{1}{n}\sum_{i\in I_n}x_i$ for all $n\in\nn$. 

Fix $n\in\nn$ and note $\|y_n\|_\xi\leqslant \frac{6}{n}$  by Lemma \ref{weaksumming}. On the other hand, $y_n$ is a convex combination of $(e_j^\xi)_{j\in \psi^{-1}(F_n)}$, thus also $Ry_n$ is a convex combination of $(e_j^\xi)_{j\in F_n}$. As $F_n\in\mathcal{S}_\xi$, it follows that $\|Ry_n\|_\xi=1$. Taking sufficiently big $n\in\nn$  we obtain the contradiction with boundedness of $R$.      
\end{proof}

\begin{rem}
The full generalization of \cite[Prop. 3.11]{GL}, i.e. $\zeta$-injectivity of $\psi$ for any $\zeta\leqslant \xi$ does not hold true in general. Take $\xi=\omega$ and the partition of $\mathbb{N}$ into successive maximal $\mathcal{S}_1$-sets $(G_i)_i$, let $m_i=\min G_i$ for each $i\in\mathbb{N}$. Let $\psi:\mathbb{N}\to M=(m_i)_i$ be defined by $\psi(e_n^\omega)=e_{m_i^\omega}$ if $n\in G_i$. Then the mapping $R:e_n^\omega\mapsto e^\omega_{\psi(n)}$, $n\in\mathbb{N}$, extends to a bounded operator, but 
$\sup_{F\in\mathcal{S}_0}\tau_0(\psi^{-1}(F))=\infty$.  However, $(e_n^\omega)_{n=1}^\infty$ dominates $(e_{m_i}^\omega)_{i=1}^\infty$ in $X_\omega$. 
\end{rem}

The next result is a version of  \cite[Prop. 3.13]{GL} for bidual Schreier spaces.

\begin{proposition}
Fix  $L,M\in[\nn]$ and $\xi<\omega_1$. Assume there is a bounded operator $V:X_\xi(L)\to Z_\xi(M)$ with 
\[\inf_{l\in L}\|Ve^\xi_l\|_0>0\]
Then there is a map $\psi:L\to M$ and a bounded operator $R:X_\xi(L)\to X_\xi(M)$ with $Re^\xi_l=e^\xi_{\psi(l)}$ for any $l\in L$.
\label{reductionoptoop2}
\end{proposition}
\begin{proof}
The reasoning follows directly the proof of \cite[Prop. 3.13]{GL} in the bidual space setting with necessary adjustments, we present it for the sake of completeness.

We recall the notion of a block diagonal matrix after \cite{LT}. Fix a matrix $(a_{ij})_{i,j=1}^\infty$. We say that a matrix $(d_{ij})_{i,j=1}^\infty$ is a block diagonal matrix of $(a_{ij})_{i,j=1}^\infty$, provided for some  $(I_n)_{n=1}^\infty, (J_n)_{n=1}^\infty\subset[\nn]^{<\infty}$ with $(I_n)_{n=1}^\infty$ and $(J_n)_{n=1}^\infty$ successive intervals we have
\[d_{ij}=\left\{\begin{array}{ll} a_{ij}, & \text{ if } (i,j)\in \bigcup_{n=1}^\infty I_n\times J_n\\ 0, & \text{ otherwise}. \end{array}\right.\] 

Let $a_{lm}=f^\xi_m(Ve^\xi_l)$ for any $l\in L$, $m\in M$, that is $Ve^\xi_l=(a_{lm})_{m\in M}\in Z_\xi(M)$ for each $l\in L$. By continuity of $V$ and the coordinate functionals $(f^\xi_i)_{i=1}^\infty$ on $Z_\xi$, for any 
$x=
\sum_{l\in L}\lambda_le^\xi_l\in X_\xi(L)$  we have \[Vx=\Big(\sum_{l\in  L}\lambda_la_{lm}\Big)_{m\in M}\in Z_\xi(M)\]
In such case we say that the matrix $(a_{lm})_{l\in L,m\in M}$ represents the operator $V\in\mathcal{L}(X_\xi(L),Z_\xi(M))$ with respect to  the basis $(e^\xi_l)_{l\in L}$ and coordinate functionals $(f^\xi_m)_{m\in M}$. 

By the assumption on $V$ there is $\delta>0$ such that for any $l\in L$ there is $m\in M$ with $|a_{lm}|>\delta$. Thus we can define a map $\psi: L\to M$ with  $|a_{l\psi(l)}|>\delta$ for any $l\in L$. As $(e^\xi_l)_{l\in L}\subset X_\xi(L)$ is weakly null, $\psi^{-1}(m)$ is finite for any $m\in M$. In particular $\psi(L)\in [M]$. 

Let $S=P_{\psi(M)}^{**}\circ V\in\mathcal{L}(X_\xi(L), Z_\xi(\psi(L)))$. Then $S$ is represented by the matrix $(a_{lm})_{l\in L, m\in \psi(L)}$ with respect to the basis $(e^\xi_l)_{l\in L}$ and coordinate functionals $(f^\xi_m)_{m\in \psi(L)}$. 
Define a matrix $(b_{lm})_{l\in L, m\in \psi(L)}$ by 
\[b_{lm}=\left\{\begin{array}{ll} a_{lm}, & \text{ if } m=\psi(l)\\ 0, & \text{ otherwise}. \end{array}\right.\] 
Notice that in each row there is a unique non-zero entry, and in each column there are finitely many non-zero entries. Therefore for some permutation $\pi$ of $L$ the matrix $(b_{\pi(l),m})_{l\in L,m\in \psi(L)}$ is a block diagonal matrix of $(a_{\pi(l),m})_{l\in L, m\in \psi(L)}$. The latter matrix represents the  operator $S\in\mathcal{L}(X_\xi(L), Z_\xi(\psi(L))$ with respect to the basis $(e_{\pi(l)})_{l\in L}$ and coordinate functionals $(f^\xi_m)_{m\in \psi(L)}$. Repeating the proof of \cite[Prop. 1.c.8]{LT} (we use  the fact that the operators on $Z_\xi$ of the form $Z_\xi\ni (\lambda_n)_{n=1}^\infty\mapsto (\pm \lambda_n)_{n=1}^\infty\in Z_\xi$ have norm 1) we obtain that the matrix  $(b_{\pi(l),m})_{l\in L, m\in M}$ defines a bounded operator $X_\xi(L)\to Z_\xi(M)$ with respect to the basis $(e_{\pi(l)})_{l\in L}$ and coordinate functionals $(f^\xi_m)_{m\in \psi(L)}$. Therefore $(b_{lm})_{l\in L,m\in \psi(L)}$ defines a bounded operator $R:X_\xi(L)\to Z_\xi(\psi(L))$ with $Re^\xi_l=e^\xi_{\psi(l)}$ for all $l\in L$. As $Re_l^\xi\in X_\xi$ for each $l\in L$,  $R(X_\xi(L))\subset X_\xi(\psi(L))\subset X_\xi(M)$, which ends the proof. 
\end{proof}

With the above preparation of tools substituting \cite[Thm 1.1]{GL}, we follow closely in the sequel the scheme of \cite[Section 4.1]{MP} in the setting of $\xi<\omega_1$. 
Fix a normalized block sequence $(u_i)_{i=1}^\infty\subset X_\xi$ 2-equivalent to the unit vector basis of $c_0$.  Let $T\in\mathcal{L}(X_\xi)$ be the operator carrying $e^\xi_i$ to $u_i$, for all $i\in\nn$.

\begin{rem}\label{remarkT2}
The operator $T^{**}\in\mathcal{L}(Z_\xi)$ satisfies  $2\|x\|_0\geq \|T^{**}x\|_\xi$, $x\in Z_\xi$. 

Indeed, note that $T=\Phi\circ id_{\xi,0}$, where $ id_{\xi,0}: X_\xi\to c_0$ is the formal identity and  $\Phi:c_0\ni (a_i)_{i=1}^\infty\mapsto \sum_{i=1}^\infty a_iu_i\in X_\xi$ is an isomorphism onto its image with $\|\Phi\|\leqslant 2$. Then $T^{**}=\Phi^{**}\circ  id_{\xi,0}^{**}$, with the formal identity $ id_{\xi,0}^{**}: Z_\xi\to\ell_\infty$ (see Prop. \ref{Zxi}) and $\Phi^{**}: \ell_\infty\to Z_\xi$ with $\|\Phi^{**}\|\leqslant 2$. Thus for any $x\in Z_\xi$ we have 
\[\|T^{**}x\|_\xi=\|\Phi^{**}( id_{\xi,0}^{**}x)\|_\xi\leqslant 2\| id_{\xi,0}^{**}x\|_0=2\|x\|_0.\]
\end{rem}

To state the next result we need dual notions to strictly and finitely strictly singular operators. Recall that an operator $S:X\to Y$, for Banach spaces $X,Y$, is strictly cosingular if for no closed infinite codimensional subspace $Z\subset Y$ the operator $S\circ \pi_Z$ is not surjective, where $\pi_Z$ is the canonical projection $\pi_Z:Y\to Y/Z$. We do not have "full" duality here; for any $S\in\mathcal{L}(X,Y)$, if $S^*$ is strictly singular (resp. cosingular), then $S$ is strictly cosingular (resp. singular), but not vice versa. The situation improves in the case of finitely strictly singular operators and finitely strictly cosingular operators (called also super strictly cosingular), introduced in \cite{P}, with the corrected definition in \cite{FHR}. Given any bounded operator $S: X\to Y$ and any closed subspace $E$ of $Y$ let
\[a_E(S)=\inf\left\{\frac{\mathrm{dist}(Sx,E)}{\mathrm{dist}(x, S^{-1}(E)}: x\in X\setminus S^{-1}(E)\right\}\]
if $S(X)\not\subset E$ and $a_E(S)=0$ otherwise. For any $n\in\nn$ let
\[ a_n(S)=\sup\{a_E(S):E\subset Y, \text{ codim}(E)=n, Y=S(X)+E \}\]
with the convention $\sup\emptyset=0$. We say that $S$ is finitely strictly cosingular provided $a_n(S)\xrightarrow{n\to\infty}0$.
By \cite[Prop. 2.13]{FHR} (see also \cite{P}) an operator $S$ is finitely strictly singular (resp. cosingular) iff its adjoint $S^*$ is finitely strictly cosingular (resp. singular). 

The next observation extends \cite[Lemma 4.14]{BKL} on the strict singularity of $T$ for the spaces $X_\xi$, $1\leqslant\xi<\omega$. 
\begin{proposition} Fix $1\leqslant\xi<\omega_1$. Then the operators $T\in\mathcal{L}(X_\xi), T^{**}\in\mathcal{L}(Z_\xi)$ are finitely strictly singular, the operator $T^*\in\mathcal{L}(X_\xi^*)$ is finitely strictly cosingular. 
\label{T finitely ss sc}
\end{proposition}
\begin{proof} By \cite[Prop. 2.13]{FHR} it is enough to prove that $T$ is finitely strictly singular. Fix $\varepsilon>0$. Pick $k\geq \varepsilon^{-1}$.  By \cite{M}, see also \cite[Lemma 13]{CP} in any $k$-dimensional subspace of $X_\xi$ there is a vector $x$ which is $k$-flat, i.e. such that $\|x\|_0=1$ and at least $k$ many coordinates equal one. As $\mathcal{S}_1\subset\mathcal{S}_\xi$ by Prop. \ref{schreierfamilyfacts} (v), for any $k$-flat $x>k$ we have $\|x\|_\xi\geq k$, whereas $\|Tx\|_\xi=\|x\|_0=1$. Therefore $n=2k$ satisfies the required condition for $\varepsilon$.     
\end{proof}

For any  $I\in[\nn]$ let $T_I=T\circ P_I\in\mathcal{FSS}(X_\xi)$. 
\begin{proposition}\label{criterion-dual}  Fix $1\leqslant \xi<\omega_1$, $I\in[\nn]$ and $\mathcal{J}\subset[\nn]$.
Assume $\|T_I^{**}-Q\|< 1$ for some $Q$ in the operator ideal generated by the family of operators $\{T_J^{**}: J\in \mathcal{J}\}$. Then there are $J_1,\dots,J_k\in\mathcal{J}$ and a  $\xi$-injective  map  $\psi: I\to J_1\cup\dots\cup  J_k$.  
\end{proposition}
\begin{proof}
We repeat the proof of \cite[Prop. 4.2]{MP} for the second adjoints of $T_I$ and $(T_J)_J$.  Assume 
$\|T_I^{**}-\sum_{t=1}^kQ_t\circ T_{J_t}^{**}\circ S_t\|<1$  for some $Q_t,S_t\in\mathcal{L}(Z_\xi)$ and $J_t\in\mathcal{J}$, $t=1,\dots,k$.  Then 
As $T_I^{**}(e^\xi_i)=e^\xi_i$, $i\in I$, the sequence $((\sum_{t=1}^kQ_t\circ T_{J_t}^{**}\circ S_t)e^\xi_i)_{i\in I}$ is seminormalized. It follows that for any $i\in I$ there is $t_i\in\{1,\dots,k\}$ so that the sequence $((T^{**}_{J_{t_i}}\circ S_{t_i})e^\xi_i)_{i\in I}$ is seminormalized. For any $t=1,\dots,k$ let $I_t=\{i\in I: t=t_i\}$ and note that $(I_t)_{t=1,\dots,k}$ form a partition of $I$. 

Let $J_0=J_1\cup\dots\cup J_k$. 
Define a bounded operator $S=\sum_{t=1}^k P^{**}_{J_0}\circ S_t\circ P_{J_t}: X_\xi(I)\to Z_\xi(J_0)$ and notice that for any $i\in I$ we have by Remark \ref{remarkT2}
\begin{align*}
\|Se^\xi_i\|_0&=\|P_{J_0}^{**}(S_{t_i}e^\xi_i)\|_0
\geq \|P^{**}_{J_{t_i}}(S_{t_i}e^\xi_i)\|_0
\geq \frac12\|T^{**}_{J_{t_i}}(S_{t_i}e^\xi_i)\|_\xi    
\end{align*}
As the sequence $((T^{**}_{J_{t_i}}\circ S_{t_i})e^\xi_i)_{i\in I}$ is seminormalized, we have $\inf_{i\in I}\|Se^\xi_i\|_0>0$. Applying Propositions \ref{reductionoptoop2} and \ref{reductionoptomap} we obtain a  $\xi$-injective  mapping $\psi: I\to J_0$, which ends the proof. 
\end{proof}

\begin{lemma}\label{sch-dyadic} Fix $\xi<\omega_1$. There is a family $(J_r)_{r\in\rr}\subset[\nn]$ such that for any pairwise different $r,r_1,\dots,r_k\in\rr$,  there is no  $\xi$-injective  map $\psi:J_r\to J_{r_1}\cup\dots\cup J_{r_k}$.
\end{lemma}
\begin{proof}
Let $\mathcal{D}$ be a dyadic tree with the lexicographic order $\leqslant_{lex}$. Let $(d_n)_{n=1}^\infty$ be the enumeration of elements of $\mathcal{D}$ according to the lexicographic order. We define inductively $(F_n)_{n=1}^\infty\subset[\nn]^{<\infty}$satisfying for any $n\in\nn$ the following:
\begin{enumerate}
    \item $\min F_{n+1}>\max F_n$,
    \item $\min F_{n+1}>\sum_{i=1}^n|F_i|$,
    \item $\tau_\xi( F_{n+1})>n\sum_{i=1}^n\tau_\xi(F_i)$,
\end{enumerate}
Given any $\mathcal{A}\subset\mathcal{D}$ let $J_\mathcal{A}=\cup_{d_n\in\mathcal{A}}F_n$. We shall prove that for any infinite $\mathcal{A},\mathcal{B}\subset\mathcal{D}$ with $\mathcal{A}\setminus\mathcal{B}$ infinite there is no $\xi$-injective map $\psi:\mathcal{A}\to\mathcal{B}$. Then taking the family $(\mathcal{J}_\mathcal{B})_{\mathcal{B} \text{ branch of }\mathcal{D} }$ will end the proof (recall that branch of a tree is a maximal linearly ordered subset of the tree). 
 
Take $\mathcal{A}=(d_{i_n})_{n\in \nn}$, $\mathcal{B}=(d_{j_n})_{n\in \nn}$ with $\mathcal{A}\setminus\mathcal{B}$ infinite, and a  map
$\psi:J_\mathcal{A}\to J_\mathcal{B}$. Fix any $N$ with $i_N\not\in (j_n)_{n\in\nn}$  and pick $M\in\nn$ with $j_M<i_N<j_{M+1}$. Let $F=\cup_{n=1}^NF_{i_n}$ and $G=\cup_{n=1}^MF_{j_n}$. Notice that by (3) and Lemma \ref{subadditivity}
\[
\tau_\xi(F)\geq \tau_\xi(F_{i_N})\geq i_N\sum_{i=1}^{j_M}\tau_\xi(F_i)\geq N\tau_\xi(G)
\]
By (2) and choice of $M$ we have $J_\mathcal{B}\setminus G>|F|$, thus
\[\tau_\xi(\psi(F)\setminus G )\leqslant \tau_1(\psi(F)\setminus G)=1\]
Thus by Lemma \ref{subadditivity} (note that $G$ is an initial segment of $J_\mathcal{B}$ by (1)) we have
\[\tau_\xi(\psi(F))\leqslant \tau_\xi(\psi(F)\setminus G)+\tau_\xi(\psi(F)\cap G)\leqslant \frac{2}{N}\tau_\xi(F)\leqslant \frac{2}{N}\tau_\xi(\psi^{-1}(\psi(F))) \]
As $N$ was arbitrarily large, by Lemma \ref{lemmaxiinjective} the map $\psi$ is not $\xi$-injective. 
\end{proof}
Now we are ready to prove the final result of this section. 
\begin{theorem}
 Fix $\xi<\omega_1$. Let $X$ denote either $X_\xi$, $X_\xi^*$ or $Z_\xi$. Then there is a family of closed operator
ideals $(I_A)_{A\subset\rr}$ on $X$, such that $I_A\subset I_B$ iff $A\subset  B$, for any $A, B \subset\rr$. In the case of $X=X_\xi$ or $Z_\xi$ the operator ideals $(I_A)_{A\subset\rr}$ can be chosen to consist of finitely strictly singular operators, in the case of $X=X_\xi^*$ - of  finitely strictly cosingular operators. 

In particular there are  $2^\mathfrak{c}$ pairwise distinct closed operator ideals,  a chain of cardinality $\mathfrak{c}$ of  closed operator ideals and
an antichain of cardinality $2^\mathfrak{c}$ of closed operator ideals on $X$. 
\label{cardinality}
\end{theorem}

\begin{proof}
Take $(J_r)_{r\in \rr}\subset [\nn]$ provided by Lemma \ref{sch-dyadic}. For any $A\subset\rr$ let 
$I_A$ be the closed operator ideal generated by $(T_{J_r})_{r\in A}$, $(T^*_{J_r})_{r\in A}$ or $(T^{**}_{J_r})_{r\in A}$ for $X=X_\xi$, $X_\xi^*$ or $Z_\xi$ respectively. Apply  Proposition \ref{criterion-dual} (in the case of $X_\xi$ or $X^*_\xi$ use the fact that $\|S\|=\|S^*\|$ for any bounded operator $S$) and Proposition \ref{T finitely ss sc} to end the proof.
\end{proof}

\end{document}